\documentclass[reqno, 10pt]{article}
\usepackage{amssymb, amsmath, amsthm, amsfonts, amscd, epsfig}
\usepackage{color}

\usepackage[english]{babel}
\usepackage{bbm}
\usepackage{graphicx, epsfig}
\usepackage{geometry,graphicx,pgfplots,caption}
\usepackage{verbatim}
\usepackage{dsfont}
\usepackage{enumerate}
\usepackage{enumitem}

\usepackage[titletoc,toc,title]{appendix}
\usepackage[normalem]{ulem}
\newtheorem{remark}{Remark}[section]
\newtheorem{assumption}{Assumption}[section]

\newtheorem{condition}{Condition}[section]
\newtheorem{theorem}{Theorem}[section]
\newtheorem{corollary}{Corollary}[section]
\newtheorem{lemma}{Lemma}[section]
\newtheorem{proposition}{Proposition}[section]
\newtheorem{definition}{Definition}[section]

\renewcommand{\d}{\,\mathrm{d}}

\numberwithin{equation}{section}
\numberwithin{figure}{section}

\renewcommand{\d}{{\rm d}}
\newcommand{\abs}[1]{\lvert {#1} \rvert}
\newcommand{\norm}[1]{\lVert {#1} \rVert}

\numberwithin{equation}{section}
\renewcommand{\d}{\,\mathrm{d}}
\allowdisplaybreaks

\title{On the consistency of the posterior distribution for nonlinear PDE parameter identification\thanks{The work of B. Jin is supported by Hong Kong RGC General Research Fund (Project 14306824),
ANR / RGC Joint Research Scheme (A-CUHK402/24), NSFC / RGC Joint Research Scheme (N\_CUHK446/25) and a start-up fund from The Chinese University of Hong Kong.}}
\date{}
\author{Yuxin Fan\thanks{Department of Mathematics, The Chinese University of Hong Kong, Shatin, N.T., Hong Kong, China (\texttt{1155205162@link.cuhk.edu.hk}, \texttt{b.jin@cuhk.edu.hk})} \and Bangti Jin\footnotemark[2]}
\date{}

\begin{document}

\maketitle

\begin{abstract}
In this work, we investigate the estimation of a parameter in PDEs using Bayesian procedures, and focus on posterior distributions constructed using  Gaussian process priors, and its variational approximation. We establish contraction rates for the posterior distribution and the variational approximation in the regime of low-regularity parameters. Specifically, the ground truth is only assumed to belong to the support space of the prior rather than to its associated reproducing kernel Hilbert space. Also we derive contraction rates for the posterior distribution constructed using the randomly truncated Gaussian process prior. The analysis relies on a delicate approximation argument that approximates low-regularity ground truths by suitable elements in the reproducing kernel Hilbert space and balances various error sources. We illustrate the general theory on three nonlinear inverse problems for PDEs.
\end{abstract}
	
\maketitle
	
\section{Introduction}

The study of inverse problems, which aims to infer unknown parameters from noisy observational data, is a cornerstone of many scientific investigations. When these parameters are functions, e.g., a spatially varying coefficient in a partial differential equation (PDE), the problem becomes infinite-dimensional, which poses significant computational and mathematical challenges. To perform inference in such settings, a wide range of rigorous methodologies has been developed, including classical regularization techniques \cite{ Engl1996Regularization,ItoJin:2015}, Bayesian approaches \cite{Stuart2010}, and more recent machine-learning-based schemes \cite{Arridge2019data-driven,HauptmannJin:2025}. The Bayesian paradigm is distinct in that it not only provides an elegant framework for uncertainty quantification via regions of high posterior probability (credible sets) but also explicitly incorporates a probabilistic model for the measurement noise, thereby more faithfully capturing real-world data acquisition processes. This approach has been widely adopted in the PDE context, since the seminal contributions \cite{Stuart2010,Dashti2017uncertainty}.

One critical theoretical aspect of Bayesian inference for such problems is the posterior contraction rate, which quantifies the speed at which the posterior distribution concentrates around the true parameter. Although theoretical foundations of posterior contraction rates in Bayesian non-parametric statistics are well-established \cite{ghosal2017,van2008rates,van2009adaptive}, the nonlinear nature of PDE parameter identifications poses substantial challenges. The general theory, as summarized in key results (e.g., \cite[Theorem 2.1]{van2008rates}), provides a high-level blueprint: the contraction at a rate $\delta_n$ is guaranteed provided that the following three key conditions are met: first, the prior must place sufficient mass on small neighborhoods of the true parameter; second,  the effective size of the parameter space is not too large in terms of its metric entropy; and third, the prior probability of the complement of a sufficiently large neighborhood of the true parameter decays exponentially fast. Verifying these conditions for concrete nonlinear inverse problems, however, can be highly nontrivial.

One breakthrough was made by Monard, Nickl and Paternain  \cite{monard2021consistent}, who demonstrated that for a class of nonlinear inverse problems, the key to fulfilling these three abstract conditions lies in the injectivity of the forward map and the availability of a stability estimate for the inverse map. These ideas were subsequently adapted by Giordano and Nickl \cite{giordano2020consistency} to derive contraction rates for the inverse conductivity problem of an elliptic PDE. They established a quantitative control for the continuity of the forward map to propagate the prior, and conduct a stability analysis for the inverse map to lift contraction rates back to the unknown parameter.
This general strategy has sparked an active line of inquiry on deriving contraction rate results for various nonlinear PDE inverse problems, e.g., the inverse scattering problem \cite{furuya2024consistency} and an inverse potential problem for the heat equation \cite{kekkonen2022consistency}. 

The above proof strategies rely on certain conditions on the forward map and the prior that are equally natural for their variational Bayesian counterparts. This theoretical framework suggests a promising path toward establishing contraction rate theories for variational methods. Indeed, general theories for the contraction of variational posteriors have begun to emerge \cite{Zhang2017ConvergenceRO}. Very recently, these ideas have been synthesized to establish variational posterior contraction rates for several nonlinear PDE inverse problems \cite{zu2024consistencyvariationalbayesianinference}.

However, a common and restrictive limitation in this line of work, when applied to PDEs with Gaussian process priors, which represents one widely studied setting, is the assumption that the true parameter lies in the reproducing kernel Hilbert space (RKHS) of the prior; see, e.g., the works \cite{10.1214/009053606000001172,Vaart2007BayesianIW,Zhang2017ConvergenceRO}. This assumption is notably strong since for Gaussian process priors with  Mat\'ern kernels, the associated RKHS is a relatively small space of smooth functions.  Crucially, as established in the fundamental work \cite{Driscoll1973}, the Gaussian process prior typically assigns zero probability to its own RKHS. This creates a significant conceptual
gap: the convergence theorems require the truth to be in a set which the prior itself deems impossible.    

In this work, we revisit the contraction rate analysis of posterior distributions resulting from Gaussian process priors under a low regularity assumption on the unknown parameter in the PDEs, in which the ground truth may have only limited Sobolev regularity. To overcome the difficulty, we approximate the ground truth by suitable elements from the RKHS and incorporate the resulting approximation error into the contraction rate analysis. We derive contraction rates for a rescaled Gaussian process prior, which are slower than that in \cite{giordano2020consistency}. Further, we extend the analysis to variational Bayesian inference, which relaxes the regularity assumption required in the existing work \cite{zu2024consistencyvariationalbayesianinference}. Additionally, we adapt the randomly truncated Gaussian prior from \cite[Section 2.2.3]{giordano2020consistency} to the low‑regularity setting. See Theorems \ref{Thm:re-scaled consis}, \ref{thm:VB-main} and \ref{thm:RT-main} for the precise statement. We illustrate the general theory on three concrete nonlinear PDE inverse problems: diffusion coefficient identification, and inverse potential problems for elliptic and subdiffusion equations in Section \ref{sec:ApplicationSection}.  We verify the necessary stability and continuity conditions, and establish posterior contraction rates under low-regularity assumptions on the true parameter. These examples demonstrate the practical feasibility of the approach in realistic low‑regularity settings.

The idea of using an approximation argument has been employed in several prior studies \cite{van2008rates,giordano2020consistency}. Using the idea, van der Vaart and van Zanten \cite{van2008rates} developed a general theory of posterior contraction rates  for classical statistical models, e.g., density estimation, classification, and regression, where the rates are governed by the position of the truth relative to the RKHS and the small-ball probabilities of the Gaussian process. Giordano and Nickl \cite{giordano2020consistency} employed the randomly truncated Gaussian process priors, which approximate the truth through finite-dimensional truncations.  Compared with the works \cite{van2008rates,giordano2020consistency}, we also adopt this strategy but for a different purpose: to derive posterior contraction rates in a low-regularity regime, and, through PDE‑related arguments, to demonstrate the feasibility of the approach.

The rest of the paper is organized as follows. In Section \ref{sec:setting}, we describe the statistical setting, the Bayesian framework with Gaussian process priors and link functions, and the variational inference approach. In Section \ref{sec:main}, we present and discuss the main results on posterior contraction rates under low-regularity assumptions. We give the proofs of the main theorems in Section \ref{sec:proofs}. In Section \ref{sec:ApplicationSection}, we present three examples to illustrate the applicability of the general theory in less regular settings.

\section{Statistical setting}\label{sec:setting}
\subsection{Statistical model and measurement scheme}

Let  $\mathcal{O}\subset \mathbb{R}^d$ be an open bounded domain with a smooth boundary $\partial\mathcal{O}$, and let $f:\mathcal{O}\to \mathbb{R}$ be a function of interest in the PDE, and $G$ be the forward map of the governing PDE, which maps the coefficient $f$ to the PDE solution. 

The PDE solution $G(f)$ is measured at a set of random design points. Specifically, let $\mu=\mathrm{d} x / |\mathcal{O}|$ (with the Lebesgue measure $\mathrm{d} x$) be the uniform distribution on $\mathcal{O}$, with $|\mathcal{O}|=\int_{\mathcal{O}} \mathrm{d} x$. The design points $\left(X_i\right)_{i=1}^N$, $N \in \mathbb{N}$, are independent and identically distributed (i.i.d.) random variables drawn from the distribution $\mu$.
At each design point $X_i$, we take a noisy measurement of $G(f)(X_i)$ and the noise is modeled by i.i.d. Gaussian random variables $\sigma W_i$, where $\sigma$ is the noise level and $W_i \overset{\mathrm{i.i.d.}}{\sim} N(0,1)$ (the standard normal distribution),  independent of the design points $X_i$'s. The observation scheme is then given by
\begin{equation}\label{eq:measurement}
    Y_i = G(f)(X_i) + \sigma W_i, \quad W_i \mathop{\sim}^{\mathrm{i.i.d.}} N(0,1), \quad i = 1, 2, \dots, N.
\end{equation}
Below we denote by $Y^{(N)} = (Y_i)_{i=1}^N$ and $X^{(N)} = (X_i)_{i=1}^N$. 

The random pairs $(Y_i, X_i)_{i=1}^N$ are i.i.d. on $\mathbb{R} \times \mathcal{O}$. We denote their common law by $P_f^i$. Then the joint law of the full dataset $D_N= (Y_i, X_i)_{i=1}^N$ defines a product probability measure
$P_f^{(N)}=\bigotimes_{i=1}^NP_f^i$, with $P_f^i=P_f^1$ for all $i$. With $\mathrm{d}y$ being the Lebesgue measure on $\mathbb{R}$, $P_f $ has the  following probability density $p_f$ with respect to the product measure $\mathrm{d}y\times \mathrm{d}\mu$:
\begin{equation}\label{modeldensity}
   p_f(y,x) = \frac{1}{\sqrt{2\pi\sigma^2}}e ^{-\frac{1}{2}\abs{y-{G}(f)(x)}^2},\quad y\in \mathbb{R}, x\in \mathcal{O},
\end{equation}
where $|\cdot|$ denotes the Euclidean norm of vectors. The goal is to infer the parameter $f$ from the data $D_N$. 

\subsection{Parameter spaces and link functions}\label{sec:Parameter spaces}
In the Bayesian approach to inverse problems \cite{Stuart2010,Dashti2017uncertainty}, we place a prior probability distribution on the unknown parameter $f$. A common choice is a Gaussian process prior, which provides a flexible prior model over functions. Gaussian processes are naturally supported on linear spaces, while the admissible parameter set $\mathcal{F}$ for the coefficient $f$ is often constrained (e.g., to positive bounded functions) and hence not linear.

Thus, we employ a bijective reparameterization of the parameter set $\mathcal{F}$ through a regular link function $\varphi$ \cite{giordano2020consistency}  so that the prior is supported on the relevant parameter space. Specifically, we express the parameter space $\mathcal{F}$ as
\begin{equation*} \mathcal{F}=\{\varphi\circ F:F \in \Theta\},
\end{equation*}
where $\Theta$ is a closed linear subspace of $L^2(\mathcal{O})$ on which we place a Gaussian process prior. 
We further introduce a normed linear subspace $\mathcal{R} \subset \Theta$, on which the prior probability measure is supported, and take $\mathcal{R} = H^{\beta}(\mathcal{O})$ for some $\beta \geq 0$ \cite[Section 2.2]{Nickl2023BayesianNS}.

Under the reparameterization, the forward model $G$ becomes 
\begin{equation}\label{eq:repara_measurement}
    Y_i = \mathcal{G}(F)(X_i) + \sigma W_i:=G(\varphi\circ F)(X_i) + \sigma W_i, \quad W_i \mathop{\sim}^{\mathrm{i.i.d.}} N(0,1), \quad i = 1, 2, \dots, N.
\end{equation}
This reformulation allows placing a Gaussian process prior on the unconstrained parameter $F\in H^\alpha(\mathcal{O})$, while ensures that the parameter $f=\varphi\circ F$ lies in the admissible parameter set $\mathcal{F}$.

\subsection{Bayesian framework: priors and posteriors}\label{sec:bayesian framework}
To complete the Bayesian model, we have to specify a prior distribution on the parameter \( F\in \Theta \).
We employ a Whittle-Mat\'{e}rn Gaussian process prior \cite[Example 25]{giordano2020consistency}. For a smoothness parameter \(\alpha > \frac{d}{2}\), the Whittle-Mat\'{e}rn process \( M = \{M(x): x \in  \mathbb R^{d}\} \) is a stationary, centered Gaussian process with the covariance kernel $K(x,x')$ given by (with $(\cdot,\cdot)_{\mathbb{R}^d}$ being Euclidean inner product)
$$
K(x,x') = \int_{\mathbb{R}^d} e^{-i(x - x', \xi)_{\mathbb{R}^d}} \mu(\d\xi), \quad \text{with} \quad \mu(\d\xi) = (1 + |\xi|^2)^{-\alpha} \d\xi.
$$

We define the prior $\Pi'$ for $F$ to be the law of the restriction of $M$ to
$\mathcal O$:
$
\Pi' = \mathcal L(M|_{\mathcal O}).
$ 
The RKHS
$\mathcal H$ of $\Pi'$ is norm–equivalent to the Sobolev space
$H^{\alpha}(\mathcal O)$ \cite[Chapter 11]{ghosal2017}.
Moreover, $M|_{\mathcal O}$ admits a version with sample paths in
$H^{\beta}(\mathcal O)$ for every $\beta < \alpha - \frac{d}{2}$
\cite[Lemma I.4]{ghosal2017}.

The posterior distribution of \(F\) given the data \(D_N =  (Y_i, X_i)_{i=1}^N\) is derived via Bayes' theorem. Since the map $(F,(y,x))\mapsto P_{F}(y,x)$ is measurable,  by Bayes' formula, the posterior distribution $\Pi'(\cdot\mid D_N)$ of $F$ arising from the data $D_N$ in the model \eqref{eq:measurement} is given by
\begin{equation}\label{eq:PostDistr}
\Pi'(B\mid D_N)= \frac{\int_B e^{\ell^{(N)}(F)}\d\Pi'(F)}
{\int_{\Theta}	e^{\ell^{(N)}(F)}\d\Pi'(F)},
\quad \forall  ~\textnormal{Borel set $B\subseteq \Theta$},
\end{equation}
where $\ell^{(N)}(F)$ is the log-likelihood function, given by
$$\ell^{(N)}(F)=-\frac{1}{2\sigma^2}\sum_{i=1}^N [Y_i-\mathcal{G}(F)(X_i)]^2.$$
The posterior distribution $\Pi'(\cdot\mid D_N)$ is the Bayesian solution to the inverse problem.

\subsection{Variational inference}
We shall also investigate variational inference. It seeks an approximate posterior distribution within a simpler, more tractable family  \(\mathcal{Q}\) of probability distributions \cite{WainwrightJordan:2008,Blei03042017}. The variational approximation \(\widehat{q} \in \mathcal{Q}\)  minimizes the Kullback-Leibler (KL) divergence to the true posterior distribution $\Pi’(\cdot\mid D_N)$:
\begin{equation}\label{eq: variational approx}
\widehat{q} = \underset{q \in \mathcal{Q}}{\operatorname{argmin}} \, D_{\text{KL}}\left( q \, \| \, \Pi’(\cdot \mid D_N) \right),
\end{equation}
where the KL divergence $D_{\text{KL}}(P_1 \| P_2)$ for two probability measures \(P_1 \ll P_2\) is defined by \cite{KullbackLeibler:1951}
$$
D_{\text{KL}}(P_1 \| P_2) = \int \log\left( \frac{{\rm d}P_1}{{\rm d}P_2} \right) {\rm d}P_1.
$$
By restricting the variational family \(\mathcal{Q}\) to a computationally convenient form,  minimizing problem \eqref{eq: variational approx} becomes more tractable, often yielding faster approximations than directly exploring the full posterior distribution $\Pi'(\cdot\mid D_N)$. Numerically, this approach has been employed in \cite{Jin:2012,GuhaWu:2015,Ray2020VariationalBayes}.

To rigorously assess the approximation quality of variational posteriors, it is useful to employ a broader class of information divergences. For any $\tau \geq 0$ and $\tau \neq 1$, the R\'{e}nyi divergence $D_\tau(P_1 \| P_2)$ between two probability measures $P_1$ and $P_2$ is defined by \cite{Renyi:1961}
$$
D_\tau(P_1 \| P_2) =
\begin{cases}
\displaystyle
\frac{1}{\tau - 1}
\log \!\int
\!\left(\frac{{\rm d}P_1}{{\rm d}P_2}\right)^{\tau - 1}
{\rm d}P_1, & \text{if } P_1 \ll P_2,\\[.5em]
+\infty, & \text{otherwise.}
\end{cases}
$$
Note that $D_\tau(P_1 \| P_2)$ converges to the KL divergence $D_{\mathrm{KL}}(P_1 \| P_2)$ as $\tau \to 1$. The R\'{e}nyi divergence provides a flexible family of information distances that generalize the KL divergence $D_{\text{KL}}(P_1 \| P_2)$ and is particularly useful for studying theoretical properties of variational posteriors. The R\'{e}nyi divergence between the induced data laws $\{ P_F : F \in \Theta \}$ is used to relate information distances on the observation space to the geometry of the forward map $\mathcal{G}(F)$ below. 

\section{Main results and discussions} \label{sec:main}
In this work, we conduct a frequentist analysis of the posterior distribution $\Pi'(\cdot\mid D_N)$, under the assumption that there exists a fixed ground truth $F_0 = \varphi^{-1} \circ f_0$ that generates the data $D_N$ in \eqref{eq:measurement}. The objective is to show that the posterior distribution $\Pi'(\cdot\mid D_N)$ contracts towards $F_0$ at a precise rate in the following sense. 

\begin{definition}\label{def:contraction rate}
A sequence $(\delta_N)_{N\in \mathbb{N}}$ is said to be a contraction rate around $F_0$ for $\Pi'(\cdot|D_N)$, for a suitable metric $d$ over $F$, if 
 $$E_{F_0}\Pi'[d(F, F_0)>\delta_N|D_N]=o(1)\quad \mbox{as }N\to \infty.$$
\end{definition}

The expectation $E_{F_0}$ is taken with respect to the data-generating distribution $P_{F_0}^{(N)}$ under the truth $F_0$. This reflects the frequentist perspective, where randomness arises solely from the sampling distribution of the data $D_N$, while $F_0$ is fixed but unknown. The definition implies that the posterior probability mass outside a $\delta_N$-neighborhood of $F_0$ is vanishingly small in the $P_{F_0}^{(N)}$-probability. We will establish contraction rates for both exact Bayesian posterior and its variational approximations.

The main challenge of nonlinear inverse problems lies in the nonlinearity of the forward map $\mathcal{G}$. To derive contraction rates, we impose the following structural conditions \cite{Nickl2023BayesianNS}. It is worth pointing out that the regularity and stability conditions on the forward map $\mathcal{G}$ require explicit estimates on growth rates of the coefficients. The Lipschitz constant $c_LM^l$ in the estimate \eqref{lip} indicates a polynomial growth of the constant with the radius $M$. The notation $(H^k(\mathcal{O}))^*$ denotes the dual space of $H^k(\mathcal{O})$.

\begin{condition}\label{condreg}
Let $\Theta \subseteq L^2(\mathcal{O})$ be a parameter space, and let $\mathcal{G}:\Theta \rightarrow L^2(\mathcal{O})$ be a measurable forward map. Suppose that there exist a normed linear subspace
$(\mathcal{R},\|\cdot\|_{\mathcal{R}})$ of $\Theta$, constants $U>0$, $ c_L>0$,  $l\geq 0$ and $k\geq 0$ such that  for every $M>1$, with
$B_{\mathcal R}(M):=\{F\in\mathcal R:\|F\|_{\mathcal R}\le M\},$ the following estimates hold:
\begin{align}\label{bound}
    &\sup_{F \in \Theta \cap B_{\mathcal{R}}(M)} \|\mathcal{G}(F)\|_{C(\overline{\mathcal{O}})} \leq U, \\ \label{lip}
    \norm{\mathcal{G}(F_1)-\mathcal{G}(F_2)}_{L^{2}(\mathcal{O})}& \leq c_LM^l \norm{F_1-F_2} _{(H^k(\mathcal{O}))^*}, \quad \forall F_1,F_2 \in B_{\mathcal{R}}(M).
      \end{align}
\end{condition}

\begin{condition}\label{condstab}
Let $\Theta$, $\mathcal{G}$ and $B_{\mathcal{R}}(M)$ be the same as Condition \ref{condreg}. Suppose that for some $\eta>0$ and all $M>1$, there exist constants $c_{\rm s}>0$ and $q\geq0$ such that
\begin{equation}\label{stab}
        \Vert F - F_0\Vert_{L^2(\mathcal{O})}\leq c_{\rm s}M^q\Vert \mathcal{G}(F)-\mathcal{G}(F_0)\Vert _{L^{2}(\mathcal{O})}^\eta, \quad F \in \Theta \cap B_{\mathcal{R}}(M).
\end{equation}
\end{condition}

\subsection{Contraction rates for Gaussian process prior}
The standard Gaussian process prior $\Pi'$ (cf. Section \ref{sec:bayesian framework}) may be insufficient to control the complexity of the nonlinear forward map $\mathcal{G}$. To address the inherent ill-posedness and nonlinearity, some additional regularization is often required. This can be achieved by an $N$-dependent `rescaling' step \cite[(2.18)]{Nickl2023BayesianNS}. Specifically, consider the rescaled prior $\Pi=\Pi_{N}$ arising from
the base prior $\Pi^\prime$: 
\begin{align}\label{eq:Prior}
\Pi_{N}=\mathcal{L}\left(F_{N}\right), \quad F_{N}=N^{-h} F',
\end{align}
where $h>0$ is some scaling parameter, and $F'\sim\Pi'$  denotes the base prior. Roughly, if the RKHS of the base prior $\Pi'$ is $\mathcal{H}$, then the rescaled prior $\Pi_N$ is a Gaussian process with an RKHS $\mathcal{H}_N$ that is $\mathcal{H}$ equipped with the strengthened norm $\|F\|_{\mathcal{H}_N} = N^{h} \|F\|_{\mathcal{H}}$ \cite[Exercise 2.6.5]{GineNickl:2015}. 
Thus, the rescaling strengthens the RKHS penalty by a factor $N^{h}$ and increases the associated quadratic penalty by a factor $N^{2h}$.
This explicit $N$-dependent regularization is analogous to Tikhonov regularization \cite{tikhonov1943stability,Engl1996Regularization,ItoJin:2015} and is crucial for taming the the ill-posedness of the inverse problem and the nonlinearity of the  map $\mathcal{G}$. The choice of the exponent $h$ will be determined by the interplay between the regularity of the prior and the smoothing property of the map $\mathcal{G}$.
 \subsubsection{Contraction rates for the posterior distribution}\label{sec:fully bayesian contraction}
Under the rescaled prior $\Pi_N$ in \eqref{eq:Prior}, we now present the first main result: the contraction rate of the exact posterior distribution $\Pi_N(\cdot\mid D_N)$ in the sense of Definition \ref{def:contraction rate}.

\begin{theorem}\label{Thm:re-scaled consis}
Suppose that the forward map $\mathcal{G} : \Theta \to L^2(\mathcal{O})$ satisfies Condition \ref{condreg} with respect to a separable normed linear subspace $(\mathcal{R}, \| \cdot \|_{\mathcal{R}})$, with the constants $U > 0$, $c_L > 0$, $l\geq 0$, and $k \geq 0$.
Let $\Pi^{\prime}$ be the base Gaussian prior on $\mathcal{R}$ with the associated RKHS $\mathcal{H}$, and let the rescaled Gaussian prior $\Pi_N$ be defined by \eqref{eq:Prior}. 
Then for $F_0 \in \mathcal{R}$, there exists $b < 0$  such that for $\delta_N = N^{b}$ and for all $r>0$, we can choose $m$ large enough such that as $N \rightarrow \infty$,
$$
P_{F_0}^N\left(\Pi_N\left(F \in \Theta,\left\|\mathcal{G}(F)-\mathcal{G}\left(F_0\right)\right\|_{L^2(\mathcal{O})} \leq m \delta_N \mid D_N\right) \leq 1-e^{-r N \delta_N^2}\right) \rightarrow 0.
$$
The admissible choice for $b$ is characterized by the system \eqref{eqn:cond-P1}--\eqref{eqn:cond-P5} in the proof; see also Proposition \ref{prop:rate}.
Moreover, if Condition \ref{condstab} holds with $\eta>0$, then
$$
P_{F_0}^N\left(\Pi_N\left(F \in \Theta \cap B_\mathcal{R}(M) ,\left\|F-F_0\right\|_{L^2(\mathcal{O})} \leq \left(m\delta_N\right)^{{\eta}} \mid D_N\right) \leq 1-e^{-r N \delta_N^2}\right)\rightarrow 0.
$$
\end{theorem}

One limitation of existing theories for the Bayesian treatment of nonlinear inverse problems is the assumption that the true parameter $F_0$ resides in the RKHS $\mathcal{H}$ of the prior \cite[Theorem 2.2.2]{Nickl2023BayesianNS}. This is a strong assumption since functions in the RKHS are smoother by a Sobolev regularity gap of order $\frac{d}{2}$ compared to  typical sample paths of the Gaussian process \cite{steinwart2019convergence}. Theorem \ref{Thm:re-scaled consis} relaxes the condition by establishing a posterior contraction rate for the ground truth $F_0$ that belongs to a larger space $\mathcal{R}$. Notably, $\mathcal{R}$ can be chosen as the Sobolev space $H^\beta(\mathcal{O})$ with $\beta$ satisfying $0 < \beta < \alpha - \frac{d}{2}$, when the RKHS of the prior is $H^\alpha(\mathcal{O})$. The technical innovation lies in constructing a sequence of approximations $F_0^{m(N)} \in \mathcal{H}$ that converge to $F_0$ in an appropriate sense, and carefully balancing the approximation error with the other terms arising in the analysis. Note that the posterior contraction rate for $\mathcal{G}(F)$ in Theorem \ref{Thm:re-scaled consis} appears suboptimal; see Remark \ref{rmk:sub-optimal} for further discussions.

\begin{remark}
Compared with the H\"{o}lder type conditional stability estimate \eqref{stab}, in practice, for many PDE parameter identifications, the following  conditional stability estimate holds
\begin{align}\label{eqn:stab-general}
    \|F-F_0\|_{L^2(\mathcal{O})}\leq \omega(\|\mathcal{G}(F)-\mathcal{G}(F_0)\|_{L^2(\mathcal{O})}),\quad \forall F_0,F \in \Theta \cap B_{\mathcal{R}}(M),
\end{align}
where the nondecreasing function $\omega:\mathbb{R}_+ \to \mathbb{R}_+$, satisfying $\lim_{s\to0^+}\omega(s)=0$, is the modulus of continuity, which may depend on ${M}$. Under this condition, there holds
$$
P_{F_0}^N\left(\Pi_N\left(F \in \Theta \cap B_\mathcal{R}(M) ,\left\|F-F_0\right\|_{L^2(\mathcal{O})} \leq \left(m\delta_N\right)^{{\eta}} \mid D_N\right) \leq 1-e^{-r N \delta_N^2}\right)\rightarrow 0.
$$
However, it is an important issue to develop an analogous theory under the stability estimate \eqref{eqn:stab-general} in the absence of the a priori constraint $F\in \Theta\cap B_\mathcal{R}(M)$.
\end{remark}

\subsubsection{Contraction rates for variational Bayes}\label{sec:variational bayesian contraction}
Now we extend the posterior contraction theory in Section \ref{sec:fully bayesian contraction} to the variational  approximation. The analysis builds on general results for posterior contraction rates \cite{Zhang2017ConvergenceRO} and adapts the framework of \cite{zu2024consistencyvariationalbayesianinference} to the present context. Let $\mathcal{Q}$ denote a suitably chosen variational family. The variational posterior is the minimizer of the KL divergence between any element of $\mathcal{Q}$ and the posterior distribution $\Pi_N(\cdot\mid D_N) $.
Theorem \ref{thm:VB-main} gives a general bound where the rate depends on the complexity of the variational family $\mathcal{Q}$ through the term $\gamma_N$. The notation $P_{F_0}^{(N)} \widehat{q}$ denotes the expected loss under the variational approximation $\widehat{q}$ and the true distribution $P_{F_0}^{(N)}$, which captures the discrepancy between the true and approximate posterior distributions.
\begin{theorem} \label{thm:VB-main}
Suppose that the forward map $\mathcal{G} : \Theta \to L^2(\mathcal{O})$ satisfies Condition \ref{condreg} with respect to a separable normed linear subspace $(\mathcal{R}, \| \cdot \|_{\mathcal{R}})$, where the constants $U > 0$, $c_L > 0$, and $k \geq 0$. 
Let $\Pi^{\prime}$ be the base Gaussian prior on $\mathcal{R}$ with the RKHS $\mathcal{H}$, and let $\Pi_N$ be the rescaled prior in \eqref{eq:Prior}. Consider the posterior distribution $\Pi_N(\cdot \mid D_N)$ based on the data $D_N = (Y_i, X_i)_{i=1}^N$ from the model \eqref{eq:measurement}. 
Suppose that $\alpha+k\geq \frac{d(l+1)}{2}$ and let $\delta_N = N^b$ for some $b < 0$ denote the posterior contraction rate of $\Pi_N(\cdot\mid D_N)$. Define
$$
\gamma_N^2 = \frac{1}{N} \inf_{q \in \mathcal{Q}} P_{F_0}^{(N)} D_{\mathrm{KL}}\big(q \,\Vert\, \Pi_N(\cdot \mid D_N)\big),
$$
and let $\widehat{q}$ denote the variational approximation defined in \eqref{eq: variational approx}. Then for any diverging sequence $M_N\to \infty$, there holds
$$
P_{F_0}^{(N)} \widehat{q}\left(\left\|\mathcal{G}(F)-\mathcal{G}\left(F_0\right)\right\|_{L^2(\mathcal{O})}^2>M_N\left(\delta_N^2+\gamma_N^2\right)\right) \rightarrow 0
.$$
Moreover, if  Condition \ref{condstab} holds for some $\eta > 0$ and $F$, $F_0\in B_\mathcal{R}(M)$ for a fixed $M>0$, then
$$ P_{F_0}^{(N)} \widehat{q}\left(\left\|F-F_0\right\|^2_{L^2(\mathcal{O})}>M_N(\delta_N^2+\gamma_N^2)^{{\eta}}\right) \rightarrow 0.$$
\end{theorem}

Theorem \ref{thm:VB-main} provides a general  contraction rate of the variational posterior $\widehat{q}$. The overall rate decomposes into two parts: the statistical error $\delta_N^2$, which is the contraction rate achievable by the true posterior $\Pi_N(\cdot\mid D_N)$, and an approximation error $\gamma_N^2$, which arises from the variational family's inability to exactly  represent the true posterior $\Pi_N(\cdot\mid D_N)$. A successful variational Bayes approach hinges on constructing a variational family $\mathcal{Q}$ that is rich enough to make the KL divergence term $\gamma_N^2$ vanish sufficiently fast yet maintains the computational efficiency, so that the approximation error $\gamma_N^2$ does not dominate the statistical error $\delta_N^2$.

One natural variational family $\mathcal{Q}$ in the  Gaussian setting is given by truncated Gaussian
measures. The motivation is threefold. First, this family is well suited to the proof strategy. The analysis relies on a localization argument around the ground truth, and truncated Gaussian measures provide a convenient way to enforce such a localization while preserving the Gaussian prior structure. Second, the choice enables direct control of the variational approximation error $\gamma_N^2$ in
Theorem \ref{thm:VB-main}. In particular,  $\gamma_N^2$ is of the same
order as the statistical error $\delta_N^2$, so that the variational approximation does not degrade
the contraction rate. Third, from a computational perspective, this family may also be viewed as a localized Gaussian approximation. Below we illustrate Theorem~\ref{thm:VB-main} using a truncated Gaussian variational family.
To this end, we define a truncated Gaussian measure $Q_N$ by conditioning the rescaled prior $\Pi_N$ on the set
$$
\widetilde{\Theta}_N
:=
\left\{
F\in\Theta:
D_{\rm KL}(P^N_{F_0}\|P^N_F)\le N\delta_N^2,\ 
\|F-F_0\|_{\mathcal R}\le M
\right\},
$$
where $M>0$ is sufficiently large. That is,
$$
Q_N(B)
=
\frac{\Pi_N(B\cap \widetilde{\Theta}_N)}{\Pi_N(\widetilde{\Theta}_N)},
\quad \forall  \mbox{ Borel set } B \subset \Theta.
$$
The next corollary shows that, if the variational family $\mathcal Q$ contains $Q_N$, then the resulting variational posterior distribution achieves the same contraction rate as the exact one.
\begin{corollary}[Truncated Gaussian priors] \label{ex:VB-truncated}
Let the conditions of Theorem \ref{thm:VB-main} hold. Let $\Pi'$ be the base prior on $\mathcal{R}$ with the associated RKHS $\mathcal{H}$, and $\Pi_N$ be the rescaled prior in \eqref{eq:Prior}. For $\delta_N=N^{b}$ with $b<0$ and the variational posterior $\widehat{q}$ defined in \eqref{eq: variational approx} with any variational set $\mathcal{Q}$ containing the truncated Gaussian measure $Q_N$, and for any diverging sequence $M_N\to \infty$, we have 
$$
P_{F_0}^{(N)} \widehat{q}\left(\left\|\mathcal{G}(F)-\mathcal{G}\left(F_0\right)\right\|_{L^2(\mathcal{O})}>M_N\delta_N\right) \rightarrow 0
.$$
Moreover, if Condition \ref{condstab} also holds with $\eta > 0$ and $F$, $F_0\in B_\mathcal{R}(M)$ for a fixed $M>0$, then 
$$ P_{F_0}^{(N)} \widehat{q}\left(\left\|F-F_0\right\|_{L^2(\mathcal{O})}>M_N\delta_N^{\eta}\right) \rightarrow 0.$$
\end{corollary}

\subsubsection{Extension of the general theory}\label{ssec:extension}

Theorems \ref{Thm:re-scaled consis} and \ref{thm:VB-main} are given under the scaled Gaussian prior $\Pi_N$. In some PDE models, the conditional stability constant $C_{\rm s}M^q$ in condition \ref{condstab} may grow too fast in the sieve radius $R_N=MN^\rho$; see Remark \ref{rmk:restriction towards q} for further discussions.

To overcome this issue, we may employ a conditional Gaussian prior
\begin{align}\label{eq:conditional prior}
\Pi_{N,M}(A):=\Pi_N(A\mid F\in\Theta_M),
\quad \mbox{with } \Theta_M=\{F\in\mathcal R:\|F\|_{\mathcal R}\le M\}, 
\end{align}
for a fixed $M>0$ chosen sufficiently large so that $F_0\in\Theta_M$. 
The boundedness assumption is frequently employed in the literature \cite{abraham2020statisticalcalderonproblems, nickl2018bernsteinvonmises}. The assumption that the parameter $F_0$ lies in a fixed bounded subset of $\mathcal{R}$ is natural and poses no essential restriction.
The prior $\Pi_{N,M}$ satisfies $\Pi_{N,M}(\Theta_M)=1$ and is absolutely continuous with respect to $\Pi_N$ with its Radon-Nikodym derivative proportional to $\mathbf 1_{\Theta_M}$. Since conditional stability is then applied on the fixed ball $\Theta_M$, the parameter contraction rate becomes $\delta_N^{1/\eta}$ without any additional blow-up. This is precisely the posterior contraction rate on $F$ stated in Theorems \ref{Thm:re-scaled consis} and \ref{thm:VB-main}.

The verification of conditions \ref{con:KL-V} to \ref{con:covering_number} in the proof of the general contraction
theorem in Section \ref{sec:fully bayesian contraction} carries over from $\Pi_N$ to $\Pi_{N,M}$. Indeed, \ref{con:KL-V} can be proved with the same contraction rate $\delta_N$, cf. \eqref{eq:star}. \ref{con:outside_prob} holds trivially since $\Pi_N(\Theta_M)$ has a positive lower bound. \ref{con:covering_number} is unaffected since the covering number bound depends only on the geometry of the sieve and not on the prior.   This reasoning applies also to the proof of Theorem \ref{thm:VB-main}. 

\subsection{Contraction rates for randomly truncated Gaussian process prior}

Now we employ a randomly truncated Gaussian process prior which involves a random dimensionality parameter $J$ with a hyper-prior with suitable tail behavior \cite[Section 2.2.3]{giordano2020consistency} instead of the Gaussian process prior rescaled in an $N$-dependent way. Let $
\{\Psi_{\ell r}:\ell\geq -1,\ r\in\mathbb Z^d\}
$
be an orthonormal basis of $L^2(\mathbb R^d)$ consisting of sufficiently regular,
compactly supported Daubechies wavelets \cite[Section 4.2.3]{GineNickl:2015}. We impose the following interior support condition on the ground truth $F_0$. Let
$\mathcal{K}\Subset \mathcal{K}'\Subset \mathcal{O}$ be two compact subsets of $\mathcal O$, and assume that
$$
F_0\in H^\beta_{\mathcal K}(\mathcal O)
:=
\{F\in H^\beta(\mathcal{O}):\operatorname{supp}(F)\subset \mathcal{K}\}.
$$
Let $\chi\in C_c^\infty(\mathcal O)$ be a cut-off function such that
$\chi=1$ on $\mathcal{K}'$.
For each $\ell\geq -1$, let
$
\mathcal{E}_\ell:=\{r\in\mathbb Z^d:
\operatorname{supp}(\Psi_{\ell r})\cap \mathcal{K} \neq\varnothing\}
$, i.e., the set of indices $r$ for which the support of $\Psi_{\ell r}$ intersects nontrivially with $\mathcal{K}$. Then there holds $|\mathcal{E}_\ell|\lesssim 2^{\ell d}$ since $\mathcal O$ is bounded \cite[Theorem 4.2.10]{GineNickl:2015}. 

For $J\geq -1$, set
\begin{equation}\label{eq:sieve_HJ}
\mathcal H_J:=
\operatorname{span}
\left\{
\chi\Psi_{\ell r}:\ell\leq J,\ r\in\mathcal E_\ell
\right\}
\subset L^2(\mathcal O).    
\end{equation}
For any fixed $J$, let $\Pi_J$ denote the conditional law on $\mathcal H_J$
of 
\begin{equation}\label{eq:RT-prior}
F^{(J)}
=\chi
\sum_{\ell\leq J,\ r\in\mathcal E_\ell}
2^{-\ell\beta}F_{\ell r}\Psi_{\ell r},
\quad
F_{\ell r}\stackrel{\mathrm{i.i.d.}}{\sim}N(0,1).
\end{equation}
The randomly truncated Gaussian prior is the mixture
$$
\Pi(\cdot)
=\sum_{j\geq -1}P(J=j)\Pi_j(\cdot),
$$
where the random truncation level $J$ is independent of the coefficients $F_{\ell r}$ and
satisfies the following properties
\begin{equation}\label{eq:RT-J-tail}
P(J>j)= e^{-c_1 2^{jd}\log(2^{jd})},\quad \forall j\geq 1, \quad \mbox{and}\quad
P(J=j)\ge e^{-c_2 2^{jd}\log(2^{jd})}
\end{equation}
for all sufficiently large $j$ and some constants $c_1,c_2>0$.
With $(\cdot,\cdot)_{\mathbb{R}^d}$ being the $L^2(\mathbb{R}^d)$ inner product, for $F_0\in H^\beta_{\mathcal{K}}(\mathcal O)$, define the finite-dimensional
approximation
\begin{equation}\label{eq:PHJ_F0}
P_{\mathcal H_J}F_0
:=
\chi
\sum_{\ell\leq J,\ r\in\mathcal R_\ell}
(F_0,\Psi_{\ell r})_{L^2(\mathbb R^d)}
\Psi_{\ell r}.
\end{equation}

Then we have the following contraction rates for a low-regularity parameter $F_0$ under the randomly truncated Gaussian process prior.

\begin{theorem}\label{thm:RT-main}
Suppose that the forward map $\mathcal{G} : \Theta \to L^2(\mathcal{O})$ satisfies Condition \ref{condreg} with respect to a separable normed linear subspace $\mathcal{R}=H^\beta(\mathcal O)$ for some $\beta\geq0$. Let $\Pi$ be the prior defined in \eqref{eq:RT-prior}--\eqref{eq:RT-J-tail}. Suppose  
$F_0\in H^\beta_{\mathcal{K}}(\mathcal O)$. Then with $\delta_N:=N^{-\frac{\beta+k}{2\beta+2k+d}}\sqrt{\log N}$, for every $r>0$, there exists $m>0$ sufficiently large such that
$$
P_{F_0}^N\!\left(
\Pi\!\left(
F\in\Theta:\|\mathcal G(F)-\mathcal G(F_0)\|_{L^2(\mathcal O)}\le m\delta_N
\,\middle|\, D_N
\right)
\le 1-e^{-rN\delta_N^2}
\right)\to 0.
$$
Moreover, if Condition \ref{condstab} holds with $\eta>0$, then
$$
P_{F_0}^N\left(\Pi\left(F \in \Theta \cap B_\mathcal{R}(M) ,\left\|F-F_0\right\|_{L^2(\mathcal{O})} \leq \left(m\delta_N\right)^{{\eta}} \mid D_N\right) \leq 1-e^{-r N \delta_N^2}\right)\rightarrow 0.
$$
\end{theorem}

Using the randomly truncated Gaussian process prior, Giordano and Nickl  \cite[Section 2.2.3]{giordano2020consistency} derived a posterior contraction rate for a nonlinear elliptic inverse problem with a relatively high-regularity parameter. Moreover, by selecting a random truncation level $J$, the resulting posterior distribution can achieve the minimax contraction rate for $\mathcal{G}(F)$ up to a logarithmic factor, thereby providing a fully Bayesian regularization mechanism that avoids an explicit dependence on the sample size $N$. Theorem \ref{thm:RT-main} extends their result to a low-regularity parameter setting. Note that unlike the scaled Gaussian process priors in Theorems \ref{Thm:re-scaled consis} and \ref{thm:VB-main}, the contraction rate for $\mathcal{G}(F)$ in Theorem \ref{thm:RT-main} actually retains the minimax polynomial order.

\section{Proofs of Theorems \ref{Thm:re-scaled consis}, \ref{thm:VB-main} and \ref{thm:RT-main}}\label{sec:proofs}
In this section, we provide the proofs of the main results in Section \ref{sec:main}.
\subsection{Preliminary results}
First, we recall several useful results. The first result \cite[Theorem 1.2]{Linde1999} characterizes the ``small ball'' asymptotics of a Gaussian measure $\mu$ on a separable Banach space $\mathcal{B}$, which are quantified by the function
$$
\phi(\varepsilon)=-\log \mu\left(x \in \mathcal{B}:\|x\|_{\mathcal{B}} \leq \varepsilon\right)
$$
for $\varepsilon>0$, in terms of the metric entropy 
$$
\mathcal{H}(\varepsilon)=\log N\left(\left\{x \in \mathcal{B}:\|x\|_{\mathcal{H}} \leq 1\right\},\|\cdot\|_{\mathcal{B}}, \varepsilon\right).
$$
Here, the notation $N\left(\left\{x \in \mathcal{B}:\|x\|_{\mathcal{H}} \leq 1\right\},\|\cdot\|_{\mathcal{B}}, \varepsilon\right)$ denotes the covering number of the unit ball in the RKHS $\mathcal{H}$ with respect to the norm $\|\cdot\|_{\mathcal{B}}$, i.e., the minimal number of $\|\cdot\|_{\mathcal{B}}$-balls of radius $\varepsilon$ required to cover $\left\{x \in \mathcal{B}:\|x\|_{\mathcal{H}} \leq 1\right\}$.

\begin{lemma}\label{lem:Small-deviation}
Let $a>0$. Then as $\varepsilon \rightarrow 0$,
$$
\phi(\varepsilon) \approx \varepsilon^{-a} \quad \text { if and only if } \quad \mathcal{H}(\varepsilon) \approx \varepsilon^{-\frac{2 a}{2+a}} \text {, }
$$
where $\approx$ denotes two-sided inequalities up to multiplicative constants (that may depend on $\mathcal{B}$). 
\end{lemma}

The next lemma \cite[Section 3.3.2]{Edmunds_Triebel_1996} gives a metric entropy bound for Sobolev balls measured in dual Sobolev norms. It is a key tool in verifying the complexity condition for posterior contraction.
\begin{lemma}\label{lem:metric-entropy}
Let $\mathcal{O}\subset\mathbb{R}^d$ be a smooth bounded domain in $\mathbb{R}^d$. Let $H^\alpha(\mathcal{O})$ denote the Sobolev space of order $\alpha$. Let $h^\alpha(r)$ denote a ball of radius $r > 0$ in $H^\alpha(\mathcal{O})$. Then for any $k > 0$, the metric entropy $\log N\left(h^\alpha(r),\|\cdot\|_{(H^k(\mathcal{O}))^*}, \epsilon\right)$ of the ball $h^\alpha(r) $ with respect to the dual norm \( \|\cdot\|_{(H^k(\mathcal{O}))^*} \) satisfies
$$
\log N\left(h^\alpha(r),\|\cdot\|_{(H^k(\mathcal{O}))^*}, \epsilon\right) \lesssim \left( r\epsilon^{-1} \right)^{\frac{d}{\alpha + k}}, \quad \forall\, 0 < \epsilon < r.
$$
\end{lemma}

The third result \cite[Theorem 2.6.12]{GineNickl:2015} gives a bound on the Gaussian measure of a fixed set, in which the associated RKHS $\mathcal{H}$ plays a crucial role. 
\begin{lemma}\label{lem:Gaussian-Isoperimetry}
Let $O_{\mathcal{H}}$ be the unit ball centred at zero of the RKHS $\mathcal{H}$ of $X$, where $X$ is a centred Gaussian random variable on a separale Banach space $\mathcal{B}$. Let $\mu$ be the probability law of $X$, a probability measure on the Borel $\sigma$-algebra $\mathfrak{A}_{\mathcal{B}}$ of $\mathcal{B}$. Then, for every set $A \in \mathfrak{A}_{\mathcal{B}}$ and every $\varepsilon>0$,
$$
\mu\left(A+\varepsilon O_{\mathcal{H}}\right) \geq \Phi\left(\Phi^{-1}(\mu(A))+\varepsilon\right),
$$
where $\Phi$ is the standard normal cumulative distribution function.
\end{lemma}

To control the behavior of the supremum of Gaussian priors,  we employ a fundamental result -- Fernique's theorem \cite[Theorem 2.1.20]{GineNickl:2015}. It provides moment estimates on the supremum and establishes an inequality relating the supremum to the bound on the process variances.
\begin{lemma}\label{lem:Fernique's theorem}
Suppose that $\{X(t), t \in T\}$ is a separable centered Gaussian process satisfying
$P\{\sup_{t \in T}|X(t)|<\infty\}$
$>0.$
Then the following two statements hold.
\begin{itemize}
\item[{\rm(a)}] $\sigma=\sigma(X):=\sup\limits_{t \in T}\left(E X^2(t)\right)^{\frac12}<\infty$ and $E \sup\limits_{t \in T}|X(t)|<\infty$;
\item[{\rm(b)}] With $K=\frac{1}{\pi^2}$, the following inequality holds 
$$
P\left\{ \left| \sup_{t \in T} |X(t)| - E \sup_{t \in T} |X(t)| \right| > u \right\} \leq 2 e^{-\frac{K u^2}{2\sigma^2}}.
$$
\end{itemize}
\end{lemma}

The Gaussian correlation inequality  decouples the joint probability of two symmetric convex events under a Gaussian measure into a product of their marginal probabilities \cite[Theorem B.1.2]{Nickl2023BayesianNS}.
\begin{lemma}\label{lem:Gaussian correlation inequality}
Let $\mu$ be a centred Gaussian measure on a separable Banach space $\mathcal{B}$ and let $C_1, C_2$ be two convex closed symmetric subsets of $\mathcal{B}$. Then
$\mu\left(C_1 \cap C_2\right) \geq \mu\left(C_1\right) \mu\left(C_2\right) .$   
\end{lemma}

The next lemma provides a key estimate on the truncated Karhunen-Lo\`eve approximation of the ground truth $F_0$. Suppose that $F_0$ admits the Karhunen--Lo\`eve expansion
$$
F_0 = \sum_{k=1}^\infty (F_0, \phi_k)\, \phi_k,
$$
where $(\phi_k)_{k \geq 1}$ are the $L^2(\mathcal{O})$ orthonormal eigenfunctions associated with the Mat\'{e}rn covariance operator of the Gaussian process prior and $(\cdot,\cdot)$ denotes the $L^2(\mathcal{O})$ inner product. Then the $m$-term truncated Karhunen-Lo\`{e}ve approximation of $F_0$ is defined by $$F_0^m = \sum_{k=1}^m (F_0, \phi_k)\, \phi_k.$$ 
Crucially, for each fixed $m$, the approximation $F_0^m$ lies in the RKHS $\mathcal{H}$ of the Gaussian process prior. The following inverse type estimate is essential for establishing posterior contraction rates when the ground truth $F_0$ has lower regularity than the RKHS $\mathcal{H}$ of the Gaussian process prior.
\begin{lemma}\label{lem:truncated}
Let $F_0^m $ be the $m$-term truncated Karhunen-Lo\`eve approximation of $F_0$. Then for $\alpha-\frac{d}{2}>\beta\geq 0$, there holds
$$\|F_0^m\|_{H^\alpha(\mathcal{O})}  \lesssim m^{\frac{\alpha-\beta}{d}} \|F_0\|_{H^\beta(\mathcal{O})}.$$
\end{lemma}
\begin{proof}
The Mat\'ern covariance operator has $L^2(\mathcal{O})$ orthonormal eigenfunctions $\{\phi_k\}$ with the eigenvalues $\lambda_k \asymp k^{-\frac{2\alpha}{d}}$ \cite[Example 6.17]{Stuart2010}. The $H^\alpha(\mathcal{O})$ norm can be equivalently expressed in terms of $(\phi_k)_{k=1}^\infty$ \cite[Theorem 4.1]{van2008reproducing}:
$$
\|F_0^m\|_{H^\alpha(\mathcal{O})}^2 = \sum_{k=1}^m \frac{|(F_0, \phi_k)|^2}{\lambda_k} \asymp \sum_{k=1}^m k^{\frac{2\alpha}{d}} |(F_0, \phi_k)|^2.
$$
Moreover, the eigenfunction characterization holds for all intermediate Sobolev orders $0<\beta<\alpha$ by interpolation of weighted $\ell^2$-spaces \cite[Corollary 3.2]{chandler2015interpolation}.
Thus for $F_0 \in H^\beta(\mathcal{O})$,
\begin{align*}
\|F_0^m\|_{H^\alpha(\mathcal O)}^2
&\asymp
\sum_{k=1}^m
k^{\frac{2(\alpha-\beta)}{d}}
k^{\frac{2\beta}{d}}|(F_0,\phi_k)|^2\le
m^{\frac{2(\alpha-\beta)}{d}}
\sum_{k=1}^m
k^{\frac{2\beta}{d}}|(F_0,\phi_k)|^2
\lesssim
m^{\frac{2(\alpha-\beta)}{d}}
\|F_0\|_{H^\beta(\mathcal O)}^2.
\end{align*}
This completes the proof of the lemma.
\end{proof}

Also we recall the notion of a KL-type neighborhood around the truth $F_0$ \cite[Definition 1.4]{castillo2024bayesiannonparametricstatisticsstflour}. It plays a crucial role in characterizing the convergence behavior of the posterior distribution, since it simultaneously controls both KL divergence and the variance term.
\begin{definition}[KL-type neighborhood]\label{def:KL_neighborhood}
Let $P_{F}$ and $P_{F_0}$ be the probability distributions induced by the
parameters $F$ and $F_0$, respectively.  
For two probability measures $P$ and $Q$ with $P \ll Q$, the variance $V(P,Q)$ of the log-likelihood ratio is defined by
$$
V(P,Q)
    = \int\!\left(
        \log\frac{\d P}{\d Q}
        - D_{\mathrm{KL}}(P\|Q)
      \right)^{2} \d P.
$$
For any sequence $\delta_N>0$, the \emph{KL-type neighborhood} of $F_0$ is
defined by
$$
B_{\mathrm{KL}}(F_0,\delta_N)
    = \Bigl\{
        F :
        D_{\mathrm{KL}}\!\left(P_{F_0}\middle\|P_{F}\right)
            \le N\delta_N^2,\;
        V\!\left(P_{F_0},P_F\right)
            \le N\delta_N^2
      \Bigr\}.
$$
\end{definition}

\subsection{Proof of Theorem \ref{Thm:re-scaled consis}}\label{ssec:Proof Bayes case}

On the space $\Theta$, we define the semi-metric
$$d_{\mathcal{G}}\left(F, F^{\prime}\right):=\left\|\mathcal{G}(F)-\mathcal{G}\left(F^{\prime}\right)\right\|_{L^2(\mathcal{O})}.$$
For any subset $\Theta^{\prime} \subset \Theta$, we denote by $N\left(\Theta^{\prime}, d_{\mathcal{G}}, \eta\right)$, $ \eta>0$, the covering number of the set $\Theta'$ with respect to the semi-metric $d_\mathcal{G}$.
We define a sequence of Borel regularization sets $\Theta_N\subset\Theta$ by
$$
\Theta_N=\left\{F \in \mathcal{R}: F=F_1+F_2,\,\left\|F_1\right\|_{(H^k(\mathcal{O}))^*} \leq M \delta_N,\,\left\|F_2\right\|_{\mathcal{H}} \leq M N^{\rho},\,\|F\|_{\mathcal{R}} \leq MN^{\rho}\right\},
$$
where the parameter $\rho\geq0$ balances the scaling in \eqref{eq:Prior} and is to be determined. 
By the general posterior contraction theorem from \cite[Theorem 1.3.2]{Nickl2023BayesianNS}, it suffices to verify the following three conditions:
\begin{enumerate}[label=(C\arabic*)]
\item\label{con:KL-V}$\Pi_N\left(\left\{F \in \Theta: \|\mathcal{G}(F)-\mathcal{G}(F_0)\|_{L^2(\mathcal{O})} \leq \delta_N,\|\mathcal{G}(F)\|_{\infty} \leq U\right\}\right) \geq e^{-A N \delta_N^2}$ for some $A>0$;
\item\label{con:outside_prob}$
\Pi_N\left(\Theta_N^c\right) \leq e^{-B N \delta_N^2}$  for some $B>A+2$;
\item$\label{con:covering_number}\log N\left(\Theta_N, d_{\mathcal{G}}, \bar{t} \delta_N\right) \leq N \delta_N^2$  for all $\bar{t}>0$ large enough.
\end{enumerate}
The proof involves balancing several sources of errors through the choice of key parameters in the three conditions. Let $m(N) = N^\zeta$ (with $\zeta>0$) be the truncation level in the Karhunen-Lo\`{e}ve expansion, and let $\delta_N = N^b$ (with $b<0$) be the posterior contraction rate. The parameter $h > 0$ is the rescaling parameter, and $\rho \geq 0$ is an auxiliary parameter.
These parameters should satisfy the following system of (in)equalities, which ensures that the necessary prior mass concentration conditions and error balances hold:
 \begin{align}
\tag{P1} \quad & b  > -\zeta \cdot \frac{k + \beta}{d}, \label{eqn:cond-P1}\\
\tag{P2} \quad & \zeta \cdot \frac{2(\alpha - \beta)}{d}  + 2h \leq 1 + 2b, \label{eqn:cond-P2}\\
\tag{P3} \quad & -(b + h)\left( \frac{2d}{2\alpha + 2k - d} \right) \leq 1 + 2b, \label{eqn:cond-P3}\\
\tag{P4} \quad & b + \frac{1}{2} = \rho + h, \label{eqn:cond-P4}\\
\tag{P5} \quad & (b - \rho)\left( -\frac{d}{\alpha + k} \right) +\frac{\rho ld}{\alpha + k}  \leq 1 + 2b.\label{eqn:cond-P5}
\end{align}
Throughout we assume the true function $F_0\in H^\beta(\mathcal{O})$ with $\beta \geq 0$ satisfying $\beta < \alpha - \frac{d}{2}$, while the RKHS $\mathcal{H}$ of the prior $\Pi'$ corresponds to the smoother space $H^\alpha(\mathcal{O})$.
An algebraic analysis confirms that a solution $(b, \zeta, \rho)$ to the system \eqref{eqn:cond-P1}--\eqref{eqn:cond-P5} with the requisite sign conditions ($b < 0$, $\zeta > 0$, $\rho \ge 0$) always exists, provided that $0 < h < \frac{1}{2}$. We divide the lengthy proof into four steps.

\medskip
\noindent\underline{Step 1}. Since $\mathcal{R}$ is separable, the Hahn-Banach theorem implies that the norm $\|\cdot\|_\mathcal{R}$ can be represented as
$$\left\|F^{\prime}\right\|_{\mathcal{R}}=\sup\limits_{T \in \mathcal{T}}\left|T\left(F^{\prime}\right)\right|, \quad \forall F^{\prime} \in \mathcal{R},
$$
where $\mathcal{T}$ is a countable family of continuous linear forms on the space $(\mathcal{R},\|\cdot\|_{\mathcal{R}})$. Then for $F^{\prime} \sim \Pi^{\prime}$, the collection of random variables $\{T(F'): T \in \mathcal{T}\}$ defines a centred Gaussian process with a countable index set and by hypothesis,
$$\operatorname{P}\Big(\left\|F^{\prime}\right\|_{\mathcal{R}}=\sup_{T \in \mathcal{T}}\left|T\left(F^{\prime}\right)\right|<\infty\Big)=1.$$ Thus Lemma \ref{lem:Fernique's theorem} implies that $ E\left\|F^{\prime}\right\|_{\mathcal{R}} \leq D$ 
for some constant $D$ that depends only on the base prior $\Pi^{\prime}$. Moreover, the same lemma yields a sub-Gaussian concentration inequality for the random variable 
$\|F'\|_{\mathcal{R}}$: there exists a constant $c>0$ such that
\begin{equation}\label{eqn:Fernique}
\Pi'\!\left(\|F'\|_{\mathcal{R}} - E\|F'\|_{\mathcal{R}} > u\right)
\le e^{-c u^2}, \quad u>0.
\end{equation}
From equality \eqref{eqn:cond-P4}, we have $ N^{\rho+h}=\sqrt{N}\delta_N$. $E\left\|F^{\prime}\right\|_{\mathcal{R}}$ is bounded, so for large $N$ we may take $u = \frac{M \cdot N^{\rho+h}}{2}$ in \eqref{eqn:Fernique} to derive %[{\color{magenta}add more details for the second inequality below?}]
\begin{align*}
\Pi\left(\|F\|_{\mathcal{R}}>MN^{\rho}\right) & =\Pi^{\prime}\left(\left\|F^{\prime}\right\|_{\mathcal{R}}>M \cdot N^{\rho+h}\right) \\
& \leq \Pi^{\prime}\left(\left\|F^{\prime}\right\|_{\mathcal{R}}-E\left\|F^{\prime}\right\|_{\mathcal{R}}>\frac{M \cdot N^{\rho+h}}{2}\right) \\
& \leq e^{-c M^2 N\delta_N^2} \rightarrow 0,  \quad \text{for }N \to +\infty.
\end{align*}

\noindent\underline{Step 2}. In this step, we address the small ball computation \cite[Theorem 2.6.13]{GineNickl:2015}, which characterizes the KL-type neighborhood of $F_0$. Since $F_0$ belongs to $\mathcal{R}$, we have for all $M>0$, 
\begin{equation*}
\|F-F_0\|_{\mathcal{R}}\leq M
\Rightarrow \|F\|_{\mathcal{R}}\leq \overline{m}:=M +\|F_0\|_{\mathcal{R}},
\end{equation*}
and by the estimate \eqref{bound} in Condition \ref{condreg}, there also holds $\|\mathcal{G}(F)\|_{C(\overline{\mathcal{O}})} \leq U$ for some $U=U_{\mathcal{G}}(\overline{m})$. Next, using \eqref{lip} in Condition \ref{condreg}, we have
\begin{align}
& \Pi\left(F:\left\|\mathcal{G}(F)-\mathcal{G}\left(F_0\right)\right\|_{L^2(\mathcal{O})} \leq \delta_N,\|\mathcal{G}(F)\|_{C(\overline{\mathcal{O}})} \leq U\right) \notag\\
\geq& \Pi\left(F:\left\|\mathcal{G}(F)-\mathcal{G}\left(F_0\right)\right\|_{L^2(\mathcal{O})} \leq \delta_N,\left\|F-F_0\right\|_{\mathcal{R}} \leq M\right) \notag \\
 \geq & \Pi(F:\left\|F-F_0\right\|_{{(H^k(\mathcal{O}))^*}} \leq 2\widetilde{\delta}_{N,\overline{m}},\left\|F-F_0\right\|_{\mathcal{R}} \leq M), \label{eq:star}
\end{align}
with $\widetilde{\delta}_{N,\overline{m}} = \frac{\delta_N}{2L_\mathcal{G}(\overline{m})}.$
For $m =  N^\zeta$, let $F_0^m := \sum_{k=1}^m (F_0, \phi_k)\phi_k$ be the $m$-term Karhunen-Lo\`eve approximation of $F_0$. For the measurable set ${\rm I}:=\{F:\|F-F_0\|_{{(H^k(\mathcal{O}))^*}}\leq 2\widetilde{\delta}_{N,\overline{m}}\}$, we have the following inclusion relation
$$
{\rm I}\supseteq 
\{F:\|F-F_0^{m(N)}\|_{(H^k(\mathcal{O}))^*} \leq \widetilde{\delta}_{N,\overline{m}}, \|F_0-F_0^{m(N)}\|_{(H^k(\mathcal{O}))^*} \leq \widetilde{\delta}_{N,\overline{m}}\}.
$$
By the spectral characterization of Sobolev norms, we have $\|F_0-F_0^{m(N)}\|_{(H^k(\mathcal{O}))^*} \lesssim m(N)^{-\frac{\beta+k}{d}}$ for sufficiently large $N$, and then by condition \eqref{eqn:cond-P1}, we have
\begin{align*}
\mathbf{1}\{\|F_0-F_0^{m(N)}\|_{(H^k(\mathcal{O}))^*} \leq \widetilde{\delta}_{N,\overline{m}}\} =1.
\end{align*} 
Hence ${\rm I}\supseteq \{F:\|F-F_0^{m(N)}\|_{(H^k(\mathcal{O}))^*} \leq \widetilde{\delta}_{N,\overline{m}}\}$.
Moreover, let $\widetilde{M}:=M-2\|F_0\|_{\mathcal{R}}$, we have
\begin{equation*}
\{F:\left\|F-F_0\right\|_{\mathcal{R}} \leq M\}\supseteq \{F:\|F-F_0^{m(N)}\|_{\mathcal{R}} \leq \widetilde{M}\}.    \end{equation*}
Now we split the probability in \eqref{eq:star} using Gaussian correlation inequality in Lemma \ref{lem:Gaussian correlation inequality}. Since the Gaussian measure of a shifted set can be bounded below by the centered measure times an exponential penalty depending on the shift  \cite[Corollary 2.6.18]{GineNickl:2015}, we derive
\begin{align*}
&\Pi(F:\left\|F-F_0\right\|_{{(H^k(\mathcal{O}))^*}} \leq 2\widetilde{\delta}_{N,\overline{m}},\left\|F-F_0\right\|_{\mathcal{R}} \leq M)\\
\geq&\Pi(F:\|F-F_0^{m(N)}\|_{{(H^k(\mathcal{O}))^*}} \leq \widetilde{\delta}_{N,\overline{m}},\|F-F_0^{m(N)}\|_{\mathcal{R}} \leq \widetilde{M})\\
\geq &e^{-\frac{\|F_0^{m(N)}\|^2_{H^\alpha(\mathcal{O})}\cdot N^{2h}}{2}}
\Pi\left(F:\|F\|_{(H^{{k}}(\mathcal{O}))^*}\leq \widetilde{\delta}_{N,\overline{m}}, \|F\|_{\mathcal{R}} \leq \widetilde{M}\right)\\
\geq&e^{-\frac{\|F_0^{m(N)}\|^2_{H^\alpha(\mathcal{O})}\cdot N^{2h}}{2}}
\Pi\left(F:\|F\|_{(H^{{k}}(\mathcal{O}))^*}\leq \widetilde{\delta}_{N,\overline{m}}\right)\cdot \Pi\left(F:\|F\|_{\mathcal{R}} \leq \widetilde{M}\right):={\rm II}
\end{align*}
The probability $\Pi\left(\|F\|_{\mathcal{R}} \leq \widetilde{M}\right)$ has a positive lower bound since $$ \Pi\left(\|F\|_{\mathcal{R}} \leq \widetilde{M}\right)= 1-\Pi\left(\|F\|_{\mathcal{R}} > \widetilde{M}\right)\geq1-e^{-c\widetilde{M}^2N^{2h}}\geq\tfrac{1}{2}.$$
Using this estimate, Lemma \ref{lem:truncated}  and inequality \eqref{eqn:cond-P2}, we obtain 
\begin{align*}
 {\rm II} \geq &\frac{1}{2}e^{-\frac{\|F_0^{m(N)}\|^2_{H^\alpha(\mathcal{O})}\cdot N^{2h}}{2}}
\Pi\left(F:\|F\|_{(H^{{k}}(\mathcal{O}))^*}\leq \widetilde{\delta}_{N,\overline{m}}\right) \\
\geq& \frac{1}{2} e^{-c_0\frac{\|F_0\|^2_{H^\beta(\mathcal{O})} \cdot m(N)^{\frac{2(\alpha-\beta)}{d}}\cdot N^{2h}}{2}}
\Pi\left(F:\|F\|_{(H^{{k}}(\mathcal{O}))^*}\leq \widetilde{\delta}_{N,\overline{m}}\right)\\
\geq& e^{-c_1 N \delta_N^2} \cdot
\Pi\left(F:\|F\|_{(H^k(\mathcal{O}))^*}\leq \widetilde{\delta}_{N,\overline{m}}\right),
\end{align*}
for some constant $c_0>0$ and sufficiently large $N$. 
By Lemma~\ref{lem:metric-entropy}, the ball \( h^\alpha(r) \) of radius \( r \) in \( \mathcal{H} \subset H^\alpha \) satisfies the following entropy bound: for any \( 0 < \epsilon < r \),
\begin{equation*}
\log N\left(h^\alpha(r),\|\cdot\|_{(H^k(\mathcal{O}))^*}, \epsilon \right) \lesssim \epsilon^{-\frac{d}{\alpha+k}}.
\end{equation*}
Hence, by the relationship between the small ball probability and the metric entropy in Lemma \ref{lem:Small-deviation} with $\frac{2a}{2+a}=\frac{d}{\alpha+k}$ (in the ambient space $(H^{k}(\mathcal{O}))^*\supset L^2(\mathcal{O})$), and  condition \eqref{eqn:cond-P3}, we have 
\begin{align*}
\Pi(F:\|F\|_{(H^k(\mathcal{O}))^*}\leq \widetilde{\delta}_{N,\overline{m}})&=\Pi'(F': \|F'\|_{(H^k(\mathcal{O}))^*}\leq \widetilde{\delta}_{N,\overline{m}}\cdot N^h)\\ &\geq e^{-c_1 (\delta_N\cdot N^h)^{-\frac{2d}{2\alpha+2k-d}}}\geq e^{-c_2 N \delta_N^2},
\end{align*}
for some constant $c_2>0.$

\noindent\underline{Step 3}. We prove condition \ref{con:outside_prob}. By Step 1, it suffices to prove
\begin{align*}
\Pi\left(F: F=F_1+F_2,\left\|F_1\right\|_{(H^k(\mathcal{O}))^*} \leq M \delta_N,\left\|F_2\right\|_{\mathcal{H}} \leq M\cdot N^{\rho}\right) \geq 1-\tfrac{1}{2} e^{-B N \delta_N^2},
\end{align*}
for $M$ large enough. Note that we can ignore the factor $\frac12$ by increasing the constant $B$. For the rescaled prior measure $\Pi_N$, it requires a lower bound on
$$
\Pi^{\prime}\left(F^{\prime}: F^{\prime}=F_1^{\prime}+F_2^{\prime},\left\|F_1^{\prime}\right\|_{(H^k(\mathcal{O}))^*} \leq M \delta_N \cdot N^h,\left\|F_2^{\prime}\right\|_{\mathcal{H}} \leq M \cdot N^{\rho+h}\right).
$$
Using Lemma \ref{lem:Small-deviation} again, we deduce that for some $c>0$,
$$
-\log \Pi^{\prime}\big(F^{\prime}:\|F^{\prime}\|_{(H^k(\mathcal{O}))^*} \leq M \delta_N \cdot N^h\big) \leq c\big(M \delta_N \cdot N^h\big)^{-\frac{2d}{2\alpha+2k-d}}.
$$
Then by taking any $M>(\frac{B}{c})^{-\frac{2\alpha+2k-d}{2d}}$ and using condition \eqref{eqn:cond-P3}, we obtain
$$
-\log \Pi^{\prime}\big(F^{\prime}:\|F^{\prime}\|_{(H^k(\mathcal{O}))^*} \leq M \delta_N \cdot N^h\big) \leq B(\delta_N \cdot N^h)^{-\frac{2d}{2\alpha+2k-d}}\leq B N \delta_N^2.
$$
Next, let $B_N=-2 \Phi^{-1}(e^{-B N \delta_N^2})$,
where $\Phi$ is the  standard normal cumulative distribution. By the inequality for $\Phi^{-1}$ \cite[Lemma K.6]{ghosal2017}, we have $\Phi^{-1}(x)\sim -\sqrt{\log(-x)}$. Thus we have $B_N \simeq 2\sqrt{B N} \delta_N$ as $N \rightarrow \infty$. Meanwhile, condition \eqref{eqn:cond-P4} ensures $N^{\rho+h}=\sqrt{N} \delta_N$. Hence for $M \geq \sqrt{B}$, there holds
\begin{align*}
& \Pi^{\prime}\left(F^{\prime}: F^{\prime}=F_1^{\prime}+F_2^{\prime},\left\|F_1^{\prime}\right\|_{(H^k(\mathcal{O}))^*} \leq M \delta_N \cdot N^h,\left\|F_2^{\prime}\right\|_{\mathcal{H}} \leq M \cdot N^{\rho+h}\right) \\
\geq &\Pi^{\prime}\left(F^{\prime}: F^{\prime}=F_1^{\prime}+F_2^{\prime},\left\|F_1^{\prime}\right\|_{(H^k(\mathcal{O}))^*} \leq M \delta_N \cdot N^h,\left\|F_2^{\prime}\right\|_{\mathcal{H}} \leq B_N\right):={\rm III}.
\end{align*}
By Lemma \ref{lem:Gaussian-Isoperimetry}, the probability ${\rm III}$ is lower bounded by
\begin{align*}
{\rm III}&\geq \Phi\left(\Phi^{-1}\left[\Pi^{\prime}\left(\left\|F^{\prime}\right\|_{(H^k(\mathcal{O}))^*} \leq M \delta_N \cdot N^h\right)\right]+B_N\right)\\
& \geq \Phi\left(\Phi^{-1}\left[e^{-B N \delta_N^2}\right]+B_N\right)
 =1-e^{-B N \delta_N^2}.
\end{align*}
\underline{Step 4}.
To prove condition \ref{con:covering_number}, we construct a $\frac{\overline{t} \delta_N }{4(MN^{\rho})^l}$-cover in the $(H^k(\mathcal{O}))^*$-distance of a ball in $\mathcal{H}$ intersected with a ball in $\mathcal{R}$, for $\overline{t}$ large enough.
Thus it suffices to prove 
$$\log N\left(h^\alpha(MN^\rho), \|\cdot\|_{(H^k(\mathcal{O}))^*},\frac{\overline{t}\delta_N}{4(MN^{\rho})^l}\right)\leq N \delta_N^2. $$
By Lemma~\ref{lem:metric-entropy} and inequality \eqref{eqn:cond-P5}, we obtain
$$\log N\left(h^\alpha(MN^\rho),\|\cdot\|_{(H^k(\mathcal{O}))^*},\frac{\overline{t}\delta_N}{4(MN^{\rho})^l}\right)\lesssim \left({N^b}{N^{-\rho}}\right)^{-\frac{d}{\alpha+k}}\cdot N^{\frac{\rho l d}{\alpha+k}}\leq N \delta_N^2.$$
Thus we have verified conditions \ref{con:KL-V}--\ref{con:covering_number}, and the desired assertion contraction rate for $\mathcal{G}(F)$ follows:
$$
P_{F_0}^N\left(\Pi_N\left(F \in \Theta,\left\|\mathcal{G}(F)-\mathcal{G}\left(F_0\right)\right\|_{L^2(\mathcal{O})} \leq m \delta_N \mid D_N\right) \leq 1-e^{-r N \delta_N^2}\right) \rightarrow 0.
$$ 
The contraction rate for $F$ follows directly from the conditional stability estimate in Condition~\ref{condstab}.

\begin{remark}
In the proof of Step 2, we lower bound the prior mass of a \((H^k(\mathcal{O}))^{*}\)-ball around \(F_{0}\). We
choose \(m(N)=N^\zeta\) so that the approximation error satisfies
$\|F_{0}-F_{0}^{m(N)}\|_{(H^k(\mathcal{O}))^{*}}
\ \le\ \frac{\delta_{N}}{2L_{\mathcal G}(\bar m)}$,
which is enforced by condition \eqref{eqn:cond-P1}. Conditions \eqref{eqn:cond-P2} and \eqref{eqn:cond-P3} are precisely the requirements that the two exponents satisfy \(\lesssim N\delta_{N}^{2}\), which are needed in the general posterior contraction theorem, while conditions \eqref{eqn:cond-P4}–\eqref{eqn:cond-P5} control the complexity of the regularization set. 
\end{remark}

\begin{remark}\label{rmk:sub-optimal}
The system of algebraic constraints \eqref{eqn:cond-P1}-\eqref{eqn:cond-P5} governing the posterior contraction rate $\delta_N = N^b$ implies a lower bound derived primarily from the intersection of constraints \eqref{eqn:cond-P1}, \eqref{eqn:cond-P2} and \eqref{eqn:cond-P5}. In PDE settings in which $l$ is negligible, we typically have the dominant restriction 
$$
b > -\frac{\beta+k}{2\alpha + 2k},
$$
under the condition $\alpha>\beta+\frac{d}{2}$.
See Proposition \ref{prop:rate} for detailed discussions on the bound on the contraction rate.
This exponent exhibits a structural gap from the standard minimax rate $N^{-\frac{\beta+k}{2\beta + 2k + d}}$ whenever the RKHS regularity $\alpha$ exceeds the sample path regularity $\beta$. Such sub‑optimality is inherent to the rescaled prior framework: the smoothness mismatch $(\alpha > \beta)$ necessarily requires deviation from the optimal scaling in order to maintain prior mass around the function $F_0 \notin H^\alpha(\mathcal{O})$; the excess smoothness $\alpha-\beta$ then deteriorates the attainable contraction rate. 
This sub-optimality is also consistent with existing lower bound results for Gaussian priors in linear models  \cite[Theorem 2]{Castillo:2008}. Castillo \cite[Theorem 2]{Castillo:2008} showed that the posterior contraction rate is determined by the small-ball term and the RKHS approximation term. When the RKHS induced by the prior is substantially smoother than the ground truth, the approximation of the truth by RKHS elements becomes costly; consequently, the RKHS approximation error may dominate the concentration rate, leading to sub-optimal contraction rates. The mechanism of the  restriction in the rescaled prior setting is similar to that in \cite{Castillo:2008}.
\end{remark}

\begin{remark}\label{rmk:restriction towards q}
Theorem \ref{Thm:re-scaled consis} is proved for the semi-metric $d_{\mathcal{G}}$ on the space $\mathcal{R}$. When transferring the rate to the parameter $F$ in the $L^2(\mathcal{O})$ metric, we may restrict attention to a fixed ball $B_{\mathcal{R}}(M)$, which leads to a contraction rate of order $\delta_N^\frac{1}{\eta}$ for $F$. Alternatively, one may employ a sieve whose radius grows as $R_N \asymp N^{\rho}$, cf. the definition of $\Theta_N$, and apply \eqref{stab} with $M \sim R_N$ to derive a contraction result over  $\mathcal{R}$. This however incurs an additional factor of $N^{\rho q}$, which necessitates a small $q$ so that $q \rho + b \eta < 0$. This restriction is especially severe in the setting of the general stability estimate \eqref{eqn:stab-general}.
\end{remark}

\subsection{Proof of Theorem \ref{thm:VB-main}}

In the proof, we adopt the parameterization in Section \ref{ssec:Proof Bayes case}. Let $\delta_N = N^{b}$ and $m(N) = N^{\zeta}$ be the contraction rate and truncation level, respectively,
with auxiliary parameters $h>0$ and $\rho\geq 0$. %[{\color{magenta}should be $\rho\geq 0$??}] 
We enforce the system of constraints \eqref{eqn:cond-P1}, \eqref{eqn:cond-P2}, \eqref{eqn:cond-P4}, \eqref{eqn:cond-P5} on the parameters $(b, c, h, \rho)$ and strengthen \eqref{eqn:cond-P3} to more precisely incorporate the entropy structure:
\begin{align}
\tag*{(P3)$'$} \quad & -(b + h)\left( \frac{2d}{2\alpha + 2k - d} \right) + \frac{2\rho ld}{2\alpha+2k-d} \leq 1 + 2b.\label{eqn:cond-P3-1}
\end{align}
We adapt the argument \cite[Theorem 3.4]{zu2024consistencyvariationalbayesianinference} to the case $F_0\ \not\in \mathcal{H}$, and verify the three conditions formulated in \cite{Zhang2017ConvergenceRO}. 
Specifically, by \cite[Theorem 2.1]{Zhang2017ConvergenceRO}, it suffices to prove that the following three conditions hold for a loss function \(L(\cdot, \cdot)\) and constants \(c, c_1, c_2, c_3 > 0\) with \(c > c_2 + c_3 + 2\):
\begin{enumerate}[label=({C\arabic*})$'$]
\item\label{c1p}For any $\delta>\delta_N$, there exist a set $\Theta_N(\delta)$ and a testing function $\phi_N$ such that
$$
P_{F_0}^{(N)} \phi_N+\sup _{\substack{F \in \Theta_N(\delta) \\ L(P_F^{(N)}, P_{F_0}^{(N)}) \geq c_1 N \delta^2}} P_F^{(N)}\left(1-\phi_N\right) \leq e^{-c N \delta^2};
$$
\item\label{c2p}For any $\delta>\delta_N$, the set $\Theta_N(\delta)$ in \ref{c1p} satisfies
$\Pi\left(\Theta_N(\delta)^c\right) \leq e^{-cN \delta^2}$; 
\item\label{c3p}For some constant $\rho>1$,
$\Pi(D_\rho(P_{F_0}^{(N)} \| P_F^{(N)}) \leq c_3 N \delta_N^2) \geq e^{-c_2 N \delta_N^2}.$
\end{enumerate}
Steps 1 to 3 below verify conditions \ref{c1p} to \ref{c3p} separately. 
Let $M > 0$ be sufficiently large and $\delta > \delta_N$. Let
\begin{equation*} \label{eq:def_RN}
R_N(\delta) := \frac{M N^{\rho}\delta}{\delta_{N}}, \quad \widetilde{\delta}_N = \frac{M\delta_N}{L_\mathcal{G}(R_N(\delta))}\quad \mbox{and} \quad \widetilde{\delta}=\frac{\delta}{L_\mathcal{G}(R_N(\delta))}. 
\end{equation*}
The quantity $R_N(\delta)$ represents a crucial radius scaling with  $N$ and $\delta$. In the proof, we define the set
\begin{equation}\label{eqn:H_N} 
H_N(\delta) = \big\{F: F = F_1 + F_2, \norm{F_1}_{(H^k(\mathcal{O}))^*}\leq \widetilde{\delta}_N, \| F_2 \|_{\mathcal{H}}\leq R_N(\delta)\big\},
\end{equation}
and the parameter set
$\Theta_N(\delta) = H_N(\delta) \cap B_{\mathcal{R}}\!\left(R_N(\delta)\right)$.

\noindent \underline{Step 1}. For condition \ref{c1p}, we employ the argument of \cite[Theorem 7.1.4]{GineNickl:2015}. Recall that the Hellinger distance $h(P,Q)$ is defined by 
$$h(P,Q) = \left(\int (\sqrt{\d P} - \sqrt{\d Q})^2\right)^{\frac12},$$
and denote by $p_F$ the density of $P_F$ (the data-generating distribution of $F$). Fix $\bar{t}>0$, and for any $j\in\mathbb{N}$, let 
$$
S_j := \{F \in \Theta_N(\delta):4j \bar{t} \delta \leq h(p_{F},p_{F_0})< 4(j+1)\bar{t}\delta\}\subseteq \Theta_N(\delta).$$
Next we bound the metric entropy of the set $\Theta_N(\delta)$. 
By \cite[Proposition 2.1]{zu2024consistencyvariationalbayesianinference}, we have
$$
h(p_{F_1},p_{F_2}) \leq \tfrac{1}{2} d_{\mathcal{G}}(F_1,F_2).
$$
By combining this bound with the Lipschitz condition \eqref{lip} of $\mathcal{G}$, we obtain
\begin{align*}
      & \log N(\Theta_N(\delta),h,j\bar{t}\delta) \leq \log N(\Theta_N(\delta),h,\bar{t}\delta)\\
            \leq &\log N(\Theta_N(\delta),d_{\mathcal{G}},2\bar{t}\delta)
            \leq \log N\big(\Theta_N(\delta),\norm{\ \cdot \ }_{(H^{k})^*},2\bar{t}\widetilde{\delta}\big).
        \end{align*}
Since $\Theta_N(\delta) \subset H_N(\delta)$, we only need to bound
$\log N\big(H_N(\delta),\norm{\ \cdot \ }_{(H^k(\mathcal{O}))^*},2\bar{t}\widetilde{\delta}\big).$
The definition \eqref{eqn:H_N} of the set $H_N(\delta)$ implies that a $\bar{t}\widetilde{\delta}$-covering in $\norm{\ \cdot \ }_{(H^k(\mathcal{O}))^*}$ of $B_{\mathcal{H}}(R_N(\delta)):=\{h\in\mathcal H:\|h\|_{\mathcal H}\le R_N(\delta)\}$ is also a $2\bar{t}\widetilde{\delta}$-covering in $\norm{\ \cdot \ }_{(H^k(\mathcal{O}))^*}$ of $H_N(\delta)$ 
for $\bar{t}$ large enough.
Thus, it suffices to bound $\log N(B_{\mathcal{H}}(R_N(\delta)),\norm{\ \cdot \ }_{(H^k(\mathcal{O}))^*},\bar{t}\widetilde{\delta}).$
Since $B_{\mathcal{H}}(R_N(\delta))\subset h^\alpha(cR_N(\delta))$, by Lemma \ref{lem:metric-entropy}, we have
\begin{align}
    \log N\big(B_{\mathcal{H}}(R_N(\delta)),\norm{\ \cdot \ }_{(H^{k}(\mathcal{O}))^*},\bar{t}\widetilde{\delta}\big) \leq &\log N\big(h^{\alpha}(cR_N(\delta)),\norm{\ \cdot \ }_{(H^{k}(\mathcal{O}))^*},\bar{t}\widetilde{\delta}\big)\nonumber\\
    \leq &c\left(\frac{N^{\rho}}{\delta_N}\right)^{\frac{d}{\alpha+{k}}}\cdot\left(\frac{N^{\rho}\delta}{\delta_N}\right)^{\frac{ld}{\alpha+{k}}}.
    \label{entropybound}
\end{align}
Condition \eqref{eqn:cond-P5} implies that $b$ and $\rho$ satisfy $(b - \rho)( -\frac{d}{\alpha + k} )+\frac{\rho ld}{\alpha+k} \leq 1 + 2b$.
Further, the condition $\alpha+k\geq \frac{d(l+1)}{2}$ in Theorem \ref{thm:VB-main} implies $\frac{ld}{\alpha+k}\leq2.$ 
Thus, for any $\delta>\delta_N$ we have
\begin{align} \label{eqn:entropy-00}
\log N(S_j,h,j\bar{t}\delta) \leq c N\delta_N^2\cdot\left(\frac{\delta}{\delta_N}\right)^{\frac{ld}{\alpha+{k}}}\leq cN\delta^2.
\end{align}
With $N(\delta) := e^{cN\delta^2}$, the estimate \eqref{eqn:entropy-00} gives  $N(S_j,h,j\bar{t}\delta)\leq N(\delta)$.
Now we choose a minimum finite set $S_j^{'}$ of points in each set $S_j$ such that every $F \in S_j$ is within Hellinger distance $j\bar{t}\delta$ of at least one of these points. By the  bound \eqref{eqn:entropy-00}, for any fixed $j$, there are at most $N(\delta)$ such points $F_{j\iota} \in S_j^{'}$, and by \cite[Corollary 7.1.3]{GineNickl:2015}, for each $F_{j\iota}$, there exists a test function $\Psi_{N,j\iota}$ such that
\begin{equation*}
P^{(N)}_{F_0}\Psi_{N,j\iota} \leq e^{-c_tNj^2\bar{t}^2\delta^2}\quad \mbox{and} \quad \mathop{\sup}_{F\in S_j, h(p_{F},p_{F_{j\iota}})<j\bar{t}\delta} P^{(N)}_{F}(1-\Psi_{N,j\iota}) \leq e^{-c_tNj^2\bar{t}^2\delta^2},
\end{equation*}
for some universal constant $c_t > 0$. Let $\Psi_N = \mathop{\max}_{j,l}\Psi_{N,j\iota}$. Then, we have
\begin{align}
&P^{(N)}_{F_0}\Psi_N \leq P^{(N)}_{F_0}\Big(\sum_j\sum_\iota\Psi_{N,j\iota}\Big) \leq \sum_j\sum_\iota e^{-c_tNj^2\bar{t}^2\delta^2} \nonumber \\
\leq& N(\delta)\sum_j  e^{-c_tNj^2\bar{t}^2\delta^2}  \label{test1}
\leq \frac{1}{1-e^{-(c_t\bar{t}^2-c)}}e^{-(c_t\bar{t}^2-c)N\delta^2} 
\leq e^{-cN\delta^2},\\
& \mathop{\sup}_{\begin{array}{c}
  F \in \Theta_N(\delta)  \\
  h(p_{F},p_{F_0}) \geq 4\bar{t}\delta
    \end{array}}
    P^{(N)}_{F}(1-\Psi_N) = \mathop{\sup}_{F \in \cup_j S_j}P^{(N)}_{F}(1-\Psi_N)\leq e^{-cN \delta^2},\label{test}
\end{align}
for any $c>0$, when $M$ is large enough. Using \cite[Proposition 2.1]{zu2024consistencyvariationalbayesianinference}, and Condition \ref{condreg},
we have the following inequality for $F \in \Theta_N(\delta)$:
\begin{align*}
  h(p_{F},p_{F_0}) &\geq c_U\| \mathcal{G}(F)-\mathcal{G}(F_0)\|_{L^{2}(\mathcal{O})},
\end{align*}
with   
$c_U = \sqrt{\frac{1 - e^{-U^2 / 2}}{2 U^2}}>0$,
and $\sup_{F\in \Theta} \Vert \mathcal{G}(F) \Vert_{C(\overline{\mathcal{O}})}\le U$. Consequently,
\begin{align*}
  \left\{ F \in \Theta_N(\delta) : h(p_{F},p_{F_0}) \geq 4\bar{t}\delta \right\}
  &\supseteq 
  \left\{ F \in \Theta_N(\delta) : 
 N \Vert \mathcal{G}(F) - \mathcal{G}(F_0) \Vert^2_{L^2(\mathcal{O})} 
  \geq \frac{16\bar{t}^2 N\delta^2}{c_U^2} \right\}.
\end{align*}
This and the estimate \eqref{test} yield
\begin{equation}
  \mathop{\sup}_{ \begin{array}{c}
 {F \in \Theta_N(\delta)}  \\
  {L(P^{(N)}_{F},P^{(N)}_{F_0}) \geq c_1 N \delta^2}
  \end{array}}
   P^{(N)}_{F}(1-\Psi_N)\leq e^{-c N \delta^2},
\end{equation}
with the loss $L(P^{(N)}_{F},P^{(N)}_{F_0}) := \frac{c_1 c_U^2}{16\bar{t}^2}\,
N\|\mathcal{G}(F)-\mathcal{G}(F_0)\|_{L^{2}(\mathcal{O})}^2$. Thus, condition \ref{c1p} holds.

\medskip
\noindent\underline{Step 2}. By the definition of $\Theta_N(\delta)$, we deduce 
\begin{align*}
\Pi_N(\Theta_N(\delta)^c) &= \Pi_N(H_N(\delta)^c \cup B_{\mathcal{R}}(R_N(\delta))^c) 
 \leq \Pi_N(H_N(\delta)^c) + \Pi_N(B_{\mathcal{R}}(R_N(\delta))^c).
\end{align*}
Thus Lemma \ref{lem:Fernique's theorem} implies that $ E\left\|F^{\prime}\right\|_{\mathcal{R}} \leq D$ 
for some constant $D$ depending only on the base prior $\Pi^{\prime}$. Moreover, the same lemma implies that  there exists a constant $c>0$ such that
$$
\Pi'\!\left(\|F'\|_{\mathcal{R}} - E\|F'\|_{\mathcal{R}} > u\right)
\le e^{ -c u^2}, \quad u>0.
$$
Condition  \eqref{eqn:cond-P4} implies $\frac{N^{\rho+h}\delta}{\delta_N}=\sqrt{N}\delta$, and taking $u = M\sqrt{N}\delta$ gives
\begin{align*}
&\Pi\left(\|F\|_{\mathcal{R}}>R_N(\delta)\right)  =\Pi'(\Vert F'\Vert _{\mathcal{R}}>M\sqrt{N}\delta) \\
\leq & \Pi^{\prime}\Big(\|F^{\prime}\|_{\mathcal{R}}-E\|F^{\prime}\|_{\mathcal{R}}>\tfrac{M\sqrt{N}\delta}{2}\Big) \leq \tfrac{1}{2}e^{-cN\delta^2}.
\end{align*}
Then it suffices to prove
\begin{align*}
   \Pi_N(H_N(\delta)) \geq 1 - e^{-BN\delta^2} \geq 1 - \tfrac{1}{2}e^{-cN\delta^2},
\end{align*}
for $B = c+2$.  By Lemmas \ref{lem:Small-deviation} and \ref{lem:metric-entropy}, together with inequality \ref{eqn:cond-P3-1}, 
we obtain for some $c>0,$
\begin{align*}
&-\log\Pi'\Big(F' : \norm{F'}_{(H^{k}(\mathcal{O}))^*} \leq \widetilde{\delta}_N N^h\Big) 
\leq c\Big(\widetilde{\delta}_N N^h \Big)^{-\frac{2d}{2\alpha+2{k}-d}}\\
\leq &c\Big(\frac{M\delta_N N^h}{c_L(MN^\rho)^l}\Big)^{-\frac{2d}{2\alpha+2{k}-d}}\cdot \Big(\frac{\delta}{\delta_N}\Big)^{\frac{2ld}{2\alpha+2{k}-d}}.
\end{align*}
By combining this estimate with the assumption $\alpha+k\geq \frac{d(l+1)}{2}$, we deduce
$$
-\log\Pi'\Big(F' : \norm{F'}_{(H^{k}(\mathcal{O}))^*} \leq  \widetilde{\delta}_N\cdot N^h\Big) \leq cc_L^{\frac{2d}{2\alpha+2{k}-d}}M^{\frac{2d(l-1)}{2\alpha+2{k}-d}}N\delta^2.
$$
Fix $M \geq c_L^{\frac{1}{1-l}}\cdot\left(\frac{B}{c}\right)^{-\frac{2\alpha+2{k} - d}{2d(1-l)}}$. Then we have
\begin{equation}\label{pbabound}
   -\log\Pi'(F' : \norm{F'}_{(H^{k}(\mathcal{O}))^*} \leq  \widetilde{\delta}_N\cdot N^h) \leq BN\delta^2.
\end{equation}
Let $B_N = -2\Phi^{-1}(e^{-BN\delta^2}).$ By \cite[Lemma K.6]{ghosal2017}, we have 
$B_N \leq 2\sqrt{2BN}\delta.$ Hence, for $M \geq 2\sqrt{2B}$, we obtain
\begin{align*}
  &\Pi'\big(F' : F' = F'_1 + F'_2, \norm{F'_1}_{(H^{k}(\mathcal{O}))^*} \leq \widetilde{\delta}_N\cdot N^h, \| F'_2 \|_{\mathcal{H}}\leq M\sqrt{N}\delta\big)\\
 \geq &\Pi'\big(F' : F' = F'_1 + F'_2, \norm{F'_1}_{(H^{k}(\mathcal{O}))^*} \leq \widetilde{\delta}_N\cdot N^h, \|F'_2 \|_{\mathcal{H}}\leq B_N\big) :={\rm I}.
\end{align*}
By the estimate \eqref{pbabound} and Lemma \ref{lem:Gaussian-Isoperimetry}, we obtain
\begin{align*}
 {\rm I}  \geq& \Phi(\Phi^{-1}[\Pi'(\norm{F'_1}_{(H^{k}(\mathcal{O}))^*} \leq  \widetilde{\delta}_N \cdot N^h)] + B_N) 
   \geq \Phi(\Phi^{-1}[e^{-BN\delta^2}] + B_N) = 1 - e^{-BN\delta^2},
\end{align*}
i.e., $\Pi_N(H_N(\delta)) \geq 1 - e^{-BN\delta^2}.$ Finally, we have
$$\Pi_N(\Theta_N(\delta)^c) \leq e^{-cN\delta^2}.$$
This proves \ref{c2p}.

\medskip
\noindent \underline{Step 3}. We verify condition \ref{c3p} for $\rho=2$. Since $P_F^{(N)}$ corresponds to $N$ i.i.d.\ observations, the R\'{e}nyi divergence tensorizes as
$D_2(P^{(N)}_{F_0}\Vert P^{(N)}_F)=N\,D_2(P_{F_0}\Vert P_F)$,
by the product structure of the likelihood. Consequently,
$$
\{F:\,D_2(P_{F_0}\Vert P_F)\le c_3\delta_N^2\}
\subseteq
\{F:\,D_2(P^{(N)}_{F_0}\Vert P^{(N)}_F)\le c_3N\delta_N^2\}.
$$ By \cite[Proposition 2.1]{zu2024consistencyvariationalbayesianinference}, we have
$$D_2\left(P_{F_1} \| P_{F_2}\right)  \leq e^{4 U^2}\left\|\mathcal{G}\left(F_1\right)-\mathcal{G}\left(F_2\right)\right\|_{L^{2}(\mathcal{O})}^2,$$
where $U$ is the uniform bound from Condition \ref{condreg}.
Using the above inequality, we have 
\begin{align*}
  &\Pi_N(F : D_2(P^{(N)}_{F_0}\Vert P^{(N)}_{F})\leq c_3N\delta_N^2)\\
 \geq& \Pi_N(F : D_2(P_{F_0}\Vert P_{F})\leq c_3\delta_N^2, \norm{F-F_0}_{\mathcal{R}}\leq M)\\
 \geq& \Pi_N(F : e^{4U^2}\norm{\mathcal{G}(F)-\mathcal{G}(F_0)}^2_{L^{2}(\mathcal{O})}\leq c_3\delta_N^2, \norm{F-F_0}_{\mathcal{R}}\leq M)\\
 \geq &\Pi_N(F : \norm{\mathcal{G}(F)-\mathcal{G}(F_0)}_{L^{2}(\mathcal{O})}\leq \sqrt{c_3}e^{-2U^2}\delta_N, \norm{F}_{\mathcal{R}}\leq \overline{m}),
\end{align*}
for $\overline{m}>0$ and $M=\overline{m} + \Vert F_0\Vert _{\mathcal{R}}$. Then
following Step 1 in the proof of Theorem~\ref{Thm:re-scaled consis}, we have
\begin{align*}
&\Pi_N(F : \norm{\mathcal{G}(F)-\mathcal{G}(F_0)}_{L^{2}(\mathcal{O})}\leq \sqrt{c_3}e^{-2U^2}\delta_N, \norm{F}_{\mathcal{R}}\leq \overline{m})\\
\geq& \Pi\!\left(F:\|F-F_0\|_{(H^k(\mathcal{O}))^*}\leq \tfrac{\sqrt{c_3}e^{-2U^2}\delta_N}{L_{\mathcal{G}}(\overline{m})}, \|F\|_{\mathcal{R}}\leq\overline{m}\right) \notag\\
\geq&  \tfrac{1}{2}e^{-kN\delta_N^2}\Pi\!\left(F:\|F\|_{(H^k(\mathcal{O}))^*}\leq\tfrac{\sqrt{c_3}e^{-2U^2}\delta_N}{2L_{\mathcal{G}}(\overline{m})}\right)
\geq \tfrac{1}{2}e^{-c_2N\delta_N^2}, \end{align*}
for some constant $c_2>0$, and hence condition \ref{c3p} also holds.

Corollary \ref{ex:VB-truncated} is direct from Theorem \ref{thm:VB-main} and the following result \cite[Theorem 3.5] {zu2024consistencyvariationalbayesianinference}.
\begin{lemma}\label{lem:variational truncated}
Under the conditions of Corollary \ref{ex:VB-truncated}, there exists a truncated Gaussian measure $q_N$ in $\mathcal{Q}$ such that $\gamma_N^2 \lesssim \delta_N^2$.
\end{lemma}

\subsection{Proof of Theorem \ref{thm:RT-main}}
Now we prove Theorem \ref{thm:RT-main}. The proof follows the general strategy of \cite[Section 2.2.3]{giordano2020consistency}, with suitable modifications to accommodate the low-regularity setting.
First we derive a suitable small ball estimate for the $L^2(\mathcal{O})$-prediction risk.
\begin{lemma}\label{lem:RT-prior-mass}
Let the assumptions of Theorem \ref{thm:RT-main} hold. Let $\mathcal{H}_J$  be the finite-dimensional subspace defined in \eqref{eq:sieve_HJ} and fix $J_N \in \mathbb{N}$ such that $2^{J_N} \simeq N^\frac{1}{2\beta+2k+d}$. Let $$\delta_N:=N^{-\frac{\beta+k}{2\beta+2k+d}}\sqrt{\log N}. $$ 
Then there exists $A>0$
such that for all sufficiently large $c>0$ and all 
$N$ large enough,
$$
\Pi(
F\in\Theta:
\|\mathcal G(F)-\mathcal G(F_0)\|_{L^2(\mathcal O)}
\le c \delta_N
)
\ge e^{-A N\delta_N^2}.
$$
\end{lemma}
\begin{proof}
For $J_N \in \mathbb{N}$ such that $ 2^{J_N} \simeq N^\frac{1}{2\beta+2k+d}$, by the definition of $\delta_N$, we have 
\begin{align}\label{eqn:JN-delta}
2^{J_N d} \log 2^{J_N d} \simeq N^{\frac{d}{2 \beta+2 k+d}} \log N=N \delta_N^2.
\end{align}
Hence, by the construction of the prior in \eqref{eq:RT-prior}--\eqref{eq:RT-J-tail}, for some $c_0>0$, there holds for all $N$ large enough
$$
P\left(J=J_N\right) \gtrsim \mathrm{e}^{-2^{J_N d}  \log 2^{J_N d}} \geqslant \mathrm{e}^{-c_0 N \delta_N^2}.
$$
Then it follows that
\begin{align*}
&\Pi(F\in\Theta:\|\mathcal G(F)-\mathcal G(F_0)\|_{L^2(\mathcal O)}\le c\delta_N
)\\
\geq & \Pi_{J_N}(F\in\Theta:\|\mathcal G(F)-\mathcal G(F_0)\|_{L^2(\mathcal O)}\le c\delta_N) P(J=J_N)\\
\geq& \Pi_{J_N}(F\in\Theta:\|\mathcal G(F)-\mathcal G(F_0)\|_{L^2(\mathcal O)}\le c\delta_N) e^{-c_0N \delta_N^2}.
\end{align*}  
Next, let 
$$
P_{\mathcal{H}_{J_N}}\left(F_0\right)= \chi\sum_{\ell \leqslant J_N, r \in \mathcal{E}_{\ell}}(F_0,\Psi_{\ell r})_{L^2(\mathbb R^d)}\Psi_{\ell r}
%\left\langle F_0, \Psi_{\ell r}\right\rangle_{L^2} \Psi_{\ell r}
$$
be the $L^2(\mathbb{R}^d)$ projection of $F_0$ onto $\mathcal{H}_{J_N}$. Since $F_0\in H^{\beta}(\mathcal{O})$, from Condition \ref{condreg} and \cite[(B6)]{giordano2020consistency}, we deduce
$$
\|\mathcal{G}\left(F_0\right)-\mathcal{G}(P_{\mathcal{H}_{J_N}}\left(F_0\right))\|_{L^2(\mathcal{O})} \lesssim\|F_0-P_{\mathcal{H}_{J_N}}\left(F_0\right)\|_{(H^k(\mathcal{O}))^*} \lesssim 2^{-J_N\left(\beta+k\right)}\le \delta_N.
$$
This and the triangle inequality imply
\begin{align*}
&\Pi_{J_N}(F \left.:\left\|\mathcal{G}(F)-\mathcal{G}\left(F_0\right)\right\|_{L^2(\mathcal{O})} \leqslant c\delta_N\right) \\
 \geqslant &\Pi_{J_N}(F:\|\mathcal{G}(F)-\mathcal{G}(P_{\mathcal{H}_{J_N}}(F_0))\|_{L^2(\mathcal{O})} \leqslant c\delta_N-\|\mathcal{G}\left(F_0\right)-\mathcal{G}(P_{\mathcal{H}_{J_N}}\left(F_0\right))\|_{L^2(\mathcal{O})}) \\
 \geqslant& \Pi_{J_N}(F:\|\mathcal{G}(F)-\mathcal{G}(P_{\mathcal{H}_{J_N}}\left(F_0\right))\|_{L^2(\mathcal{O})} \leqslant c_1 \delta_N):={\rm I},
\end{align*}
for some  $c_1>0$. By \cite[(B5)]{giordano2020consistency}, we have \begin{equation}\label{eqn:P-F0} \|P_{\mathcal{H}_{J_N}}\left(F_0\right)\|_{\mathcal{H}_{J_N}} \leqslant \sqrt{\sum_{\ell \leqslant j, r \in \mathcal{E}_{\ell}} 2^{2 \ell \beta}(F_0,\Psi_{\ell r})_{L^2(\mathbb R^d)}^2} \leqslant\left\|F_0\right\|_{H^\beta(\mathcal{O})}.
\end{equation}
Hence, by Condition \ref{condreg}, we can lower bound the probability ${\rm I}$ by 
\begin{align*}
{\rm I}\geq & \Pi_{J_N}\left(F :\|\mathcal{G}(F)-\mathcal{G}(P_{\mathcal{H}_{J_N}}(F_0))\|_{L^2(\mathcal{O})} \leqslant c_1 \delta_N,\|F-P_{\mathcal{H}_{J_N}}(F_0)\|_{H^\beta(\mathcal{O})} \leqslant \delta_N\right) \\
 \geqslant & \Pi_{J_N}(F:\|F-P_{\mathcal{H}_{J_N}}(F_0)\|_{(H^k(\mathcal{O}))^*} \leqslant c_2 \delta_N,\|F-P_{\mathcal{H}_{J_N}}(F_0)\|_{H^\beta(\mathcal{O})} \leqslant \delta_N).
 \end{align*}
By the continuous embedding $\|F-P_{\mathcal{H}_{J_N}}(F_0)\|_{(H^k(\mathcal O))^*}\le \|F-P_{\mathcal{H}_{J_N}}(F_0)\|_{H^\beta(\mathcal O)}$ and choosing $c_3\le \min\{1,c_2\}$, we deduce that the $H^\beta(\mathcal{O})$-ball is contained in the $(H^k(\mathcal{O}))^*$-ball. Then by the estimate \eqref{eqn:P-F0}, we arrive at 
 \begin{align*}
{\rm I} \geq &\Pi_{J_N}(F:\|F-P_{\mathcal{H}_{J_N}}(F_0)\|_{H^\beta(\mathcal{O})} \leqslant c_3\delta_N)\\
 \geqslant& \mathrm{e}^{-\frac{1}{2}\|P_{\mathcal{H}_{J_N}}(F_0)\|_{\mathcal{H}_{J_N}}^2} \Pi_{J_N}(F:\|F\|_{H^\beta(\mathcal{O})} \leqslant c_3 \delta_N)\\
 \geqslant &\mathrm{e}^{-
\frac12\|F_0\|_{H^\beta(\mathcal{O})}^2} \Pi_{J_N}\left(F:\|F\|_{H^\beta(\mathcal{O})} \leqslant c_3 \delta_N\right),
\end{align*}
Moreover, by the wavelet characterization of the $H^{\beta}(\mathcal{O})$ norm \cite[p.~370]{GineNickl:2015}, we deduce 
$$
\begin{aligned}
 \Pi_{J_N}\left(F:\|F\|_{H^\beta(\mathcal{O})} \leqslant c_3\delta_N\right)&=P\left(\Big\| \chi\sum_{\ell \leqslant J_N, r \in \mathcal{E}_{\ell}}2^{-\ell \beta} F_{\ell r}\Psi_{\ell r}\Big\|_{H^\beta(\mathcal{O})}\leq c_3 \delta_N\right)\\
& \geq P\left(\sum_{m=1}^{{\rm dim}(\mathcal{H}_{J_N})}Z_m^2\le c_4 \delta_N^2\right) :={\rm II},
\end{aligned}
$$
for some $c_4>0$, where $Z_m\stackrel{\mathrm{i.i.d.}}{\sim}N(0,1)$. Using the inequality $P\left(\left|Z_1\right| \leqslant t\right) \gtrsim t \text { for } t \rightarrow 0$, the relation \eqref{eqn:JN-delta} and the following bound on $\operatorname{dim}(\mathcal{H}_{J_N})$:
\begin{equation}\label{eqn:dim-H}
\operatorname{dim}\left(\mathcal{H}_{J_N}\right)=\sum\limits_{l\leq J_N}|\mathcal{E}_\ell|\lesssim 2^{J_N d}\lesssim N \delta_N^2,
\end{equation}
we can lower bound the probability ${\rm II}$ by
\begin{align*}
{\rm II} \geq &P\left(\sqrt{\operatorname{dim}\left(\mathcal{H}_{J_N}\right)} \max _{m \leqslant \operatorname{dim}(\mathcal{H}_{J_N})}\left|Z_m\right| \leqslant \sqrt{c_5} \delta_N\right) \\
\geq& P\left(\max _{m \leqslant \mathrm{dim}\left(\mathcal{H}_{J_N}\right)}\left|Z_m\right| \leqslant c_5  \delta_N \operatorname{dim}\left(\mathcal{H}_{J_N}\right)^{-\frac{1}{2}}\right) \\
=& \prod_{m \leqslant \operatorname{dim}\left(\mathcal{H}_{J_N}\right)} P\left(\left|Z_m\right| \leqslant c_5 \delta_N \operatorname{dim}\left(\mathcal{H}_{J_N}\right)^{-\frac{1}{2}}\right)\\
\geq& \prod_{m \leqslant \operatorname{dim}\left(\mathcal{H}_{J_N}\right)} c_5 \delta_N \operatorname{dim}\left(\mathcal{H}_{J_N}\right)^{-\frac{1}{2}}\\
= & e^{\operatorname{dim}\left(\mathcal{H}_{J_N}\right)\operatorname{log}(c_5N^{-\frac{\beta+k+\frac{d}{2}}{2\beta+2k+d}}\sqrt{\operatorname{log}N})}\\
\geq & e^{-c_62^{J_Nd}\operatorname{log}N}=e^{-c_7N \delta_N^2}.
\end{align*}
This completes the proof of the lemma.
\end{proof}

\begin{lemma}\label{lem:RT-sieve}
For the finite-dimensional subspaces $\mathcal{H}_{J_N}$ defined in \eqref{eq:sieve_HJ}, define for each $M>0$,
$$
\mathcal A_N := \{F\in \mathcal{H}_{J_N}: \|F\|_{H^{\beta}(\mathcal{O})}\le M 2^{J_N\beta}\sqrt{N} \delta_N\}.$$
Then, for every $B>0$, there exists $M>0$ large enough such that
\begin{equation}\label{eq:sieve_out_prob}
\Pi(\mathcal A_N^c)\le e^{-B N\delta_N^2}.    
\end{equation}
Moreover, for each fixed $M>0$ and all $N$ large enough, there holds
\begin{equation}\label{eq:sieve_covering}
\log N(\delta_N;\mathcal A_N,d_{\mathcal G})
\le cN\delta_N^2.
\end{equation}
\end{lemma}
\begin{proof}
Let $Z_m\stackrel{\mathrm{i.i.d.}}{\sim} N(0,1)$. Note that $$\|F\|^2_{H^{\beta}(\mathcal{O})}\leq  \sum_{\ell \leqslant J_N, r \in \mathcal{E}_{\ell}} F_{\ell r}^2 ,\quad  \forall F \in \mathcal{H}_{J_N}.$$ Thus for any sufficiently large $N$, by the construction of the prior and \cite[Theorem 3.1.9]{GineNickl:2015} and the relations \eqref{eqn:JN-delta} and \eqref{eqn:dim-H}, we have
\begin{align*}
\Pi\left(\mathcal{A}_N^c\right) & \leqslant P\left(J>J_N\right)+P\left(\sum_{\ell \leqslant J \wedge J_N, r \in \mathcal{E}_{\ell}} F_{\ell r}^2 > M N \delta_N^2\right) \\
& \leqslant \mathrm{e}^{-2^{J_N d} \log 2^{J_N d}}+P\left(\sum_{m \leqslant \operatorname{dim}\left(\mathcal{H}_{J_N}\right)} Z_m^2>M N \delta_N^2\right) \\
&\leqslant \mathrm{e}^{-B N \delta_N^2}+P\left(\sum_{m \leqslant \operatorname{dim}\left(\mathcal{H}_{J_N}\right)}\left(Z_m^2-1\right)>\bar{M} N \delta_N^2\right)\\
&\leq  \mathrm{e}^{-B N \delta_N^2}+e^{-\frac{\bar{M}^2 (N \delta_N^2)^2}{4\operatorname{dim}(\mathcal{H}_{J_N})+4\bar{M} N \delta_N^2}}\leqslant \mathrm{e}^{-\bar{B} N \delta_N^2}.
\end{align*}
Next we derive  the estimate \eqref{eq:sieve_covering}.
Note that on the set $A_N$, the local Lipschitz constant is bounded. Hence for any $F_1,F_2\in \mathcal{A}_N$, there holds
$$
\begin{aligned}
\left\|\mathcal{G}\left(F_1\right)-\mathcal{G}\left(F_2\right)\right\|_{L^2(\mathcal{O})} & \lesssim\left(N^{\frac{\beta}{2 \beta+2k+d}} \sqrt{N} \delta_N\right)^l\|F_1-F_2\|_{(H^k(\mathcal{O}))^*} \\
& \lesssim N^{\frac{\beta l}{2\beta+2 k+d}}\left(\sqrt{N} \delta_N\right)^l \sqrt{\sum_{\ell \leqslant J_{N, r \in \mathcal{R}_l}}\left((F_1,\Psi_{\ell r})_{L^2(\mathbb{R}^d)}-(F_2,\Psi_{\ell r})_{\mathbb{R}^d}\right)^2} .
\end{aligned}
$$
Meanwhile, for any $F\in \mathcal{A}_N$, there holds \cite[(B2)]{giordano2020consistency}, $$\|F\|^2_{H^{\beta}(\mathcal{O})}\leq 2^{2 J_N\beta} \sum_{\ell \leqslant J_N, r \in \mathcal{E}_{\ell}} (F_0,\Psi_{\ell r})_{L^2(\mathbb R^d)}^2.$$ By the standard metric entropy estimate for balls $B_{\mathbb{R}^p}(r), r>0$, in the Euclidean space $\mathbb{R}^p$ \cite[Proposition 4.3.34]{GineNickl:2015} and the bound \eqref{eqn:dim-H}, we deduce that for $N$ large enough, there holds
$$
\begin{aligned}
\log N\left(\delta_N ; \mathcal{A}_N, d_{\mathcal G}\right) & \lesssim \log N\left(\delta_N N^{\frac{-\beta l}{2 \beta+2 k+d}}\left(\sqrt{N} \delta_N\right)^{-l} ; B_{\mathbb{R}^{\operatorname{dim}(\mathcal{H}_{J_N})}}\left(M \sqrt{N} \delta_N\right),\|\cdot\|_{\mathbb{R}^{\operatorname{dim}(\mathcal{H}_{J_N})}}\right) \\
& \leqslant \operatorname{dim}\left(\mathcal{H}_{J_N}\right) \log \frac{3 M \sqrt{N} \delta_N}{\delta_N N^{\frac{-\beta l}{2 \beta+2k+d}}\left(\sqrt{N} \delta_N\right)^{-l}} \lesssim N \delta_N^2 .
\end{aligned}
$$
This completes the proof of the lemma.
\end{proof}

Now we can prove Theorem \ref{thm:RT-main}.
\begin{proof}[Proof of Theorem \ref{thm:RT-main}]
By verifying the assumptions of \cite[Theorem 13]{giordano2020consistency} via Lemmas \ref{lem:RT-prior-mass} and \ref{lem:RT-sieve}, we derive the following contraction result for $\mathcal{G}(F)$:
$$
P_{F_0}^N\!\left(
\Pi\!\left(
F\in\Theta:\|\mathcal G(F)-\mathcal G(F_0)\|_{L^2(\mathcal O)}\le m\delta_N
\,\middle|\, D_N
\right)
\le 1-e^{-rN\delta_N^2}
\right)\to 0.
$$
Then using Condition \ref{condstab}, we further have
$$
P_{F_0}^N\left(\Pi\left(F \in \Theta \cap B_\mathcal{R}(M) ,\left\|F-F_0\right\|_{L^2(\mathcal{O})} \leq \left(m\delta_N\right)^{{\eta}} \mid D_N\right) \leq 1-e^{-r N \delta_N^2}\right)\rightarrow 0.
$$
This completes the proof of the theorem.
\end{proof}

\section{Applications to nonlinear PDE inverse problems}\label{sec:ApplicationSection}
In this section, we illustrate the general contraction theory on three examples (see, e.g., \cite{nickl2018bernsteinvonmises,giordano2020consistency,kekkonen2022consistency,Nickl2023BayesianNS,Siebel:2025}). 
In these inverse problems, the unknown parameter $f$ possesses only limited regularity, as often adopted in variational regularization \cite{Engl1996Regularization,ItoJin:2015}. Given a smooth bounded domain $\mathcal{O}$ and constants $0<\lambda<\Lambda<\infty$, for $\beta\geq 1$, consider the following parameter space for $f$:
$$
\mathcal{F}=\Big\{f \in L^{\infty}(\mathcal{O})\cap H^{\beta}(\mathcal{O}): \lambda<\operatorname*{ess\,inf}_{x \in \mathcal{O}} f(x)\leq \operatorname*{ess\,sup}_{x \in \mathcal{O}} f(x)<\Lambda \Big\}.
$$
To enforce the box constraints, let $\varphi:\mathbb{R}\to (\lambda,\Lambda)$ be a smooth, strictly increasing bijection with bounded derivatives of all orders, e.g., 
$\varphi(z) = \lambda + (\Lambda-\lambda)\frac{1}{1+e^{-z}}$, for $z\in\mathbb{R}$.
It parameterizes $f$ as $f = \varphi\circ F$ with an unconstrained one $F\in H^\beta(\mathcal{O})$.
For the Bayesian treatment, following the discussion in Section \ref{ssec:extension}, we fix a sufficiently large $M>0$ and consider the Sobolev ball
$
\Theta_M := \left\{F \in H^\beta(\mathcal{O}) : \|F\|_{H^\beta(\mathcal{O})} \leq M \right\}
$
such that $F_0\in \Theta_M$, and below we set 
$\Theta:=\Theta_M$. 
Throughout, the notation $(\cdot,\cdot)$ denotes the $L^2(\mathcal{O})$ inner product.

\subsection{Diffusion coefficient identification}\label{ssec:diffusion}
Let $\mathcal{O}\subset \mathbb{R}^2$ be a smooth domain and the source $g \in L^{\infty}(\mathcal{O})$ be strictly positive. In this example we take $\beta=1+\varepsilon$, with $\varepsilon\in(0,1)$, which can be arbitrarily small. Hence $f=\varphi \circ F\in H^{1+\varepsilon}(\mathcal{O})$ \cite[Lemma 2]{brezis2001composition}. 
Consider the following Darcy problem
\begin{align}\label{eqn:darcy}
\left\{\begin{aligned}-\nabla \cdot(f \nabla u)&=g,\quad \text {on } \mathcal{O}, \\ u&=0,\quad \text {on } \partial \mathcal{O} .\end{aligned}\right.
\end{align}

For $f\in \mathcal{F} $, the associated elliptic operator is uniformly elliptic, and by De Giorgi--Nash--Moser regularity theory and boundary regularity \cite[Corollary 8.28]{gilbarg1977elliptic}, for each $F\in\Theta$, the weak solution $u=\mathcal G(F)$
belongs to $C(\overline{\mathcal O})$. Moreover, by the standard elliptic regularity theory and the maximum principle, there
exists a constant $U<\infty$, depending only on
$\lambda$, $\Lambda$, $\mathcal{O}$ and $\|g\|_{L^\infty(\mathcal{O})}$, such that
$\sup_{F \in \Theta}\,\|\mathcal{G}(F)\|_{C(\overline{\mathcal{O}})}
\le U.$
The map $\mathcal{G}$ is Lipschitz continuous in the following sense.

\begin{lemma}\label{lem:diffusion_continuity}
Suppose that for some $\varepsilon\in (0,1)$,
$
\sup_{f\in\mathcal F}\|f\|_{H^{1+\varepsilon}(\mathcal O)}\leq M.
$
Then there exist $s\in(0,\varepsilon)$ and $L>0$, depending only on $\lambda$, $\Lambda$, $\mathcal O$, $\varepsilon$, $\|g\|_{L^\infty(\mathcal O)}$
and $M$, such that    
$$
\|\mathcal G(F_1)-\mathcal G(F_2)\|_{L^2(\mathcal O)}
\le L\,\|F_1-F_2\|_{(H^{s}(\mathcal O))^*},\quad \forall F_1,F_2\in \Theta.
$$
\end{lemma}
\begin{proof}
Let $u_i=\mathcal G(F_i)$ and $w=u_1-u_2$. Then $w\in H^1_0(\mathcal O)$ solves
$$
-\nabla\cdot(f_1\nabla w)=\nabla\cdot((f_2-f_1)\nabla u_2),
\quad\text{in }\mathcal O.
$$
Fix any $v\in L^2(\mathcal{O})$ and let $z_v\in H_0^1(\mathcal{O})$ solve the adjoint problem
$$
-\nabla\cdot(f_1\nabla z_v) = v,\quad \text{in }\mathcal{O}.
$$
Lax-Milgram theorem gives a unique solution $z_v\in H^1_0(\mathcal{O})$ and moreover, by Lemma \ref{lem:difference_quotient}, we have 
\begin{align}\label{eqn:reg-zv}\|\nabla z_v\|_{H^s(\mathcal{O})} \le c_2 \|v\|_{L^2(\mathcal{O})},
\end{align}
with $c_2=c_2(\lambda,\mathcal{O},M)$. By testing the equation for $w$ with $z_v$ and integration by parts, we have
$$
(w,v)
=((f_1-f_2)\nabla u_2,\nabla z_v).
$$
Since $\sup_{f\in\mathcal F}\|f\|_{H^{1+\varepsilon}(\mathcal O)}\leq M$, by Sobolev embedding \cite{Adams:2003}, 
$\sup_{f\in\mathcal F}\|f\|_{C^{0,\gamma}(\overline{\mathcal O})}\le c$ for $0<\gamma<\epsilon$. By \cite[Theorem 8.33]{gilbarg1977elliptic}, we have
$$
\|u_2\|_{C^{1,\gamma}(\overline{\mathcal O})}\le c\|g\|_{L^\infty(\mathcal{O})},
$$
where $c=c(\lambda,\Lambda,\mathcal{O},\gamma,M)$ is uniform over $\mathcal{F}$. For $s<\gamma$, Lemma \ref{multiplier estimate} and the estimate \eqref{eqn:reg-zv} yield
$$
\|\nabla u_2\cdot \nabla z_v\|_{H^s(\mathcal O)}
\le c \|\nabla u_2\|_{C^{0,\gamma}(\overline{\mathcal O})}\|\nabla z_v\|_{H^s(\mathcal O)}
\le c\|v\|_{L^2(\mathcal O)}.$$
Consequently,
$$
|(w,v)|
\le \|f_2-f_1\|_{(H^{s}(\mathcal O))^*}
   \|\nabla u_2\cdot \nabla z_v\|_{H^s(\mathcal O)}
\le c \|f_2-f_1\|_{(H^{s}(\mathcal O))^*}\|v\|_{L^2(\mathcal O)}.
$$
Taking the supremum over all $v\in L^2(\mathcal O)$ with
$\|v\|_{L^2(\mathcal O)}\le 1$ yields
$$\|w\|_{L^2(\mathcal O)}
\le c\|f_2-f_1\|_{(H^{s}(\mathcal O))^*}.
$$
Meanwhile, by the fundamental theorem of calculus and the relation $F_i=\varphi^{-1}(f_i)$, we have 
\begin{align} \label{eqn:F-fdc}
f_2-f_1=\varphi(F_2)-\varphi(F_1)=\int_0^1 \varphi'(F^\tau)\d \tau (F_2-F_1),\quad \mbox{with } F^\tau = F_1+\tau(F_2-F_1).
\end{align}
Hence, for any $\psi \in H^s(\mathcal O)$, we can write
\begin{align*}
(f_2 - f_1, \psi)
= \int_0^1 ( F_2 - F_1, \varphi'(F^\tau) \psi ) \d\tau.
\end{align*}
Therefore, by the duality pairing between $(H^{s}(\mathcal O))^*$ and $H^s(\mathcal O)$,
\begin{align*}
\|f_2 - f_1\|_{(H^{s}(\mathcal O))^*}
\le \|F_2 - F_1\|_{(H^{s}(\mathcal O))*}
\int_0^1 \sup_{\|\psi\|_{H^s(\mathcal O)} \le 1}
\|\varphi'(F^\tau) \psi\|_{H^s(\mathcal O)}  \d\tau.
\end{align*}
Next, we bound the product term $\varphi'(F^\tau)\psi$. Since $\|f_i\|_{C^{0,\gamma}(\overline{\mathcal O})}\le c$ and $\lambda \le f_i \le \Lambda$, and the function $\varphi^{-1}$ is Lipschitz in the compact interval $[\lambda,\Lambda]$, we obtain for any $x,y\in\mathcal O$,  
$$|F_i(x)-F_i(y)| = |\varphi^{-1}(f_i(x))-\varphi^{-1}(f_i(y))| \le c|f_i(x)-f_i(y)|
\le c\|f_i\|_{C^{0,\gamma}(\overline{\mathcal O})}|x-y|^\gamma.$$
Moreover, $\|F_i\|_{L^\infty(\mathcal O)}\le |\varphi^{-1}(\lambda)| + c(\|f_i\|_{L^\infty(\mathcal O)}+\lambda)$. Hence $
F_i=\varphi^{-1}(f_i)\in C^{0,\gamma}(\overline{\mathcal O}),
$ with
$$\|F_i\|_{C^{0,\gamma}(\overline{\mathcal O})}
\le |\varphi^{-1}(\lambda)| + c\lambda+c\|f_i\|_{C^{0,\gamma}(\overline{\mathcal O})} \le c,$$
and  
$\|F_\tau\|_{C^{0,\gamma}(\overline{\mathcal O})}\le c$. Since $\varphi''$ is bounded, we have for any $x,y\in\mathcal{O}$,
$$|\varphi'(F_\tau(x))-\varphi'(F_\tau(y))|
\le \|\varphi''\|_{L^\infty(\mathbb{R})} |F_\tau(x)-F_\tau(y)|
\le \|\varphi''\|_{L^\infty(\mathbb{R})} \|F_\tau\|_{C^{0,\gamma}(\overline{\mathcal O})} |x-y|^\gamma.$$
Moreover, since $\varphi'$ is globally Lipschitz,
$$
\|\varphi'(F_\tau)\|_{L^\infty(\mathcal O)}
\le |\varphi'(0)| + \|\varphi''\|_{L^\infty(\mathbb{R})} \|F_\tau\|_{L^\infty(\overline{\mathcal O})}
\le |\varphi'(0)| + \|\varphi''\|_{L^\infty(\mathbb{R})} \|F_\tau\|_{C^{0,\gamma}(\overline{\mathcal O})}.
$$
In sum, we obtain  
$$\|\varphi'(F_\tau)\|_{C^{0,\gamma}(\overline{\mathcal O})}\leq \|\varphi''\|_{L^\infty(\mathbb{R})} \|F_\tau\|_{C^{0,\gamma}(\overline{\mathcal O})}+|\varphi'(0)|\leq c.$$
Therefore, for all $\psi\in H^s(\mathcal O)$, Lemma \ref{multiplier estimate} gives
$$
\|\varphi'(F_\tau)\psi\|_{H^s(\mathcal O)}
\le c \|\varphi'(F_\tau)\|_{C^{0,\gamma}(\overline{\mathcal O})} \|\psi\|_{H^s(\mathcal O)}
\le c \|\psi\|_{H^s(\mathcal O)}.
$$
Hence $\|f_2 - f_1\|_{(H^{s}(\mathcal O))^*}
\le c \|F_2 - F_1\|_{(H^{s}(\mathcal O))^*}.$  This completes the proof of the lemma.
\end{proof}

\begin{remark}\label{rem:diffusion_dimension}
Lemma \ref{lem:diffusion_continuity} is stated for $d=2$, and the result can be extended to any dimension $d \ge 2$, for which the condition $f \in H^{1+\varepsilon}(\mathcal{O})$ should be replaced by $f \in H^{\frac{d}{2} + \varepsilon}(\mathcal{O})$ in order to ensure the H\"{o}lder continuity $f \in C^{0, \gamma}(\overline{\mathcal{O}})$ for $\gamma<\varepsilon$ via the Sobolev embedding $H^{\frac{d}{2}+\varepsilon}(\mathcal{O}) \hookrightarrow C^{0, \gamma}(\overline{\mathcal{O}})$. The analytical framework, including the use of difference quotients and the boundary flattening via diffeomorphisms, remains valid for any smooth bounded domain $\mathcal{O} \subset \mathbb{R}^d$.
\end{remark}

Hence, the map $\mathcal{G}$ satisfies Condition \ref{condreg} with the choice $k=s$ with $s>0$ arbitrarily small. 
The stability results in \cite[Corollary 3.8 and Lemma 4.1]{bonito2017diffusion} yield that for all $M>0$ and $F\in H^{1+\varepsilon}(\mathcal{O})$ with $\|F\|_{H^{1+\varepsilon}(\mathcal{O})}\le M$, there holds
$$
\|F-F_0\|_{L^2(\mathcal{O})} \lesssim \|\mathcal{G}(F)-\mathcal{G}(F_0)\|_{H^1(\mathcal{O})}^{\frac{1}{6}},
$$
which together with the Gagliardo-Nirenberg interpolation inequality \cite{BrezisMironescu:2018} implies 
\begin{align} 
\|F-F_0\|_{L^2(\mathcal{O})}&\lesssim \|\mathcal{G}(F)-\mathcal{G}(F_0)\|_{L^2(\mathcal{O})}^{\frac{s}{6(1+s)}}\|\mathcal{G}(F)-\mathcal{G}(F_0)\|_{H^{1+s}(\mathcal{O})}^{\frac{1}{6(1+s)}}\nonumber\\
&\lesssim \|\mathcal{G}(F)-\mathcal{G}(F_0)\|_{L^2(\mathcal{O})}^{\frac{s}{6(1+s)}},\label{eqn:stab-conductivity}
\end{align}
since for all $f\in \mathcal{F}$, Lemma \ref{lem:difference_quotient} implies the uniform bound
$\|\mathcal{G}(F)-\mathcal{G}(F_0)\|_{H^{1+s}(\mathcal{O})}\le c\|g\|_{L^2(\mathcal{O})}$. This shows Condition \ref{condstab}. Now we take the parameter space $\Theta\subset H^{1+\varepsilon}(\mathcal{O})$.
Let $\Pi_N$ be the Gaussian prior in Section~\ref{sec:bayesian framework}, specifically a Whittle-Mat\'{e}rn field with the RKHS $\mathcal{H}$ contained in $H^3(\mathcal{O})$ (i.e., $\alpha=3$). 
The sample paths lie in $H^{1+\varepsilon}(\mathcal{O})$ (i.e., $\beta=1+\varepsilon$).
We define a prior supported on $\Theta_M$ by
$$
\Pi_{N,M}(A):=\Pi_N(A\mid F\in\Theta_M)
=\frac{\Pi_N(A\cap\Theta_M)}{\Pi_N(\Theta_M)},\quad A\subset H^{1+\varepsilon}(\mathcal{O}).
$$
We employ $\Pi_{N,M}$ as the prior on $F$ and denote by $\Pi_{N,M}(\cdot\mid D_N)$ 
the resulting posterior distribution. By construction, $\Pi_{N,M}(\Theta_M)=1$ and 
$\Pi_{N,M}(\Theta_M^c\mid D_N)=0$ almost surely. 
For $\varepsilon\in (0,1)$, the choice 
$$\left(\beta, b, h, \zeta, \rho,k\right) = 
\left(1+\varepsilon, -\frac{2+3\varepsilon}{8(2+\varepsilon)}, \frac{2+3\varepsilon}{8(2+\varepsilon)},\frac{1}{2(2+\varepsilon)}, \frac{2-\varepsilon}{4(2+\varepsilon)},\frac{3\varepsilon}{4}\right)$$
satisfies conditions \eqref{eqn:cond-P1}--\eqref{eqn:cond-P5}, and also the strengthened condition \ref{eqn:cond-P3-1}. Thus both Theorems~\ref{Thm:re-scaled consis} and \ref{thm:VB-main} apply, and the posterior and variational posterior contract around the forward solution $\mathcal G(F_0)$ at the rate $N^{-\frac{2+3\varepsilon}{8(2+\varepsilon)}}$ in $L^2(\mathcal O)$. Since the exponent $\frac{2+3\varepsilon}{8(2+\varepsilon)}$ is increasing in $\varepsilon$, the contraction rate improves as the regularity index $\varepsilon$ increases. Under the conditional stability estimate \eqref{eqn:stab-conductivity} with the exponent $\eta=\frac{k}{6(1+k)}=\frac{\varepsilon}{6(2+\varepsilon)}$, we obtain a contraction rate for $F$.

\begin{remark}
For the model \eqref{eqn:darcy}, Giordano and Nickl \cite{giordano2020consistency} investigated the posterior consistency based on a class of Gaussian process priors, including rescaled Gaussian priors and Gaussian sieve priors, and established the Lipschitz continuity and stability estimate in a relatively regular parameter space to facilitate the PDE analysis, i.e., the conductivity $f\in H^\alpha(\mathcal O)$ with $\alpha>1+\frac{d}{2}$. Under the rescaled prior, they established a contraction rate of $N^{-\frac{\alpha+1}{2(\alpha+1)+d}}$ for $\mathcal{G}(F_0)$.
In \cite{doi:10.1137/18M1236137}, the authors constructed the MAP estimator and derived its convergence rate under the white noise setting. More recently, Siebel \cite{Siebel:2025} extended the MAP analysis to a random design regression model and proved the existence and consistency for the estimator under the assumption $f\in H^\alpha(\mathcal O)$ with $\alpha>1+\frac{5d}{2}$. 
Further, these works assume that the source term $g$ is smooth, while in this work we only require $g\in L^{\infty}(\mathcal{O})$. Thus our analysis covers a genuinely low‑regularity regime for both the unknown parameter $f$ and the known source $g$.
\end{remark}
\subsection{Inverse potential problem}\label{ssec:potential}
Let $\mathcal{O}\subset \mathbb{R}^d$  for $d\le 3$ be a smooth domain with a boundary $\partial\mathcal{O}$, and $g \in L^2(\mathcal{O})$ be strictly positive, and fix $\beta=1$. For $f \in H^1(\mathcal{O})$, consider the following problem:
$$
\left\{\begin{aligned}
		-\Delta u +fu &= g, \ &\mbox{in}&\ \mathcal{O}, \\
		u&=0, \ &\mbox{on}&\ \partial\mathcal{O}.
\end{aligned}\right.
$$
We aim at recovering the potential $f$ from the observation of $u_f$ in $\mathcal{O}$. Since $g \in L^2(\mathcal{O})$, by the standard elliptic regularity theory, $u \in H^2(\mathcal{O}) \cap H^1_0(\mathcal{O})$ and the $H^2(\mathcal{O})$-norm of 
$u$ is uniformly bounded over $\Theta$. By Sobolev embedding \cite{Adams:2003}, $u=\mathcal{G}(F) \in C(\overline{\mathcal{O}})$. Thus, there exists a constant $U < \infty$, depending only on $\lambda$, $\Lambda$, $\mathcal{O}$ and $\|g\|_{L^2(\mathcal{O})}$, such that
$\sup_{F \in \Theta}\,\|\mathcal{G}(F)\|_{C(\overline{\mathcal{O}})}
\le U.$ The map $\mathcal{G}$ is Lipschitz continuous in the following sense.
\begin{lemma}\label{ex2: Lip}
There exists a constant $L>0$, depending on $\lambda$, $\Lambda$, $\mathcal{O}$, $g$ and $\varphi$, such that
%for all $F_1,F_2\in \Theta_M$ with $f_i:=\varphi\circ F_i\in\mathcal{F}$, $i=1,2$, we have
\begin{equation}\label{eq:ex2Lip}
\bigl\|\mathcal{G}(F_1)-\mathcal{G}(F_2)\bigr\|_{L^{2}(\mathcal{O})}
\;\leq\;
L(\|F_1\|_{H^1(\mathcal{O})}+\|F_2\|_{H^1(\mathcal{O})})\|F_2-F_1\|_{H^{-1}(\mathcal{O})},\quad \forall F_1,F_2\in \Theta.
\end{equation}
\end{lemma}
\begin{proof}
Let $u_i = \mathcal{G}(F_i)$, $i=1,2$. By the standard elliptic regularity theory, we have 
$\|u_2\|_{H^2(\mathcal{O})}\leq c\|g\|_{L^2(\mathcal{O})}$,
with $c=c(\lambda, \Lambda,\mathcal{O})$. The difference $w := u_1 - u_2\in H_0^1(\mathcal{O})$ satisfies  
$$-\Delta w + f_1 w = (f_2 - f_1) u_2, \quad\text{in }\mathcal{O}.$$  
Fix any $v \in L^2(\mathcal{O})$ and let $z_v \in H_0^1(\mathcal{O})$ solve the adjoint problem
$$
-\Delta z_v + f_1 z_v = v, \quad \text{in }\mathcal{O}.
$$
Since $f_1 \geq \lambda > 0$ and $f_1 \in L^\infty(\mathcal{O})$, the elliptic regularity theory \cite[Section 6.3.2]{evans2010partial} implies the existence of a unique solution $z_v \in H^2(\mathcal{O}) \cap H_0^1(\mathcal{O})$ with $
\|z_v\|_{H^2(\mathcal{O})} \leq C \|v\|_{L^2(\mathcal{O})}$,
with $C = C(\lambda, \Lambda, \mathcal{O}) > 0$. For $d=1,2,3$, we have the following bilinear estimate \cite[Theorem 4.39]{Adams:2003}
\begin{align}
\|u_2z_v\|_{H^{2}(\mathcal{O})} \leq c\|u_2\|_{H^2(\mathcal{O})}\|z_v\|_{{H}^2(\mathcal{O})} \leq c_1\|v\|_{L^2(\mathcal{O})},\label{eqn:bilin}
\end{align}
with $c_1=c_1(\lambda,\Lambda,\mathcal{O},\|g\|_{L^2(\mathcal{O})})$.
Thus, by Sobolev embedding \cite{Adams:2003}, $u_2z_v\in H_0^1(\mathcal{O})\cap H^2(\mathcal{O}) \hookrightarrow L^\infty(\mathcal{O})$.
Testing the equations for $w$ and $z_v$ with $z_v$ and $w$, respectively, gives
$$( w,v)= (f_2 - f_1, u_2 z_v).$$
Next we appeal to the identity \eqref{eqn:F-fdc}.
The condition on the link function $\varphi$ implies $\varphi'(F^\tau)\in L^\infty(\mathcal{O})$ and $\nabla \varphi'(F^\tau) = \varphi''(F^\tau)\nabla F^\tau\in L^2(\mathcal{O})$.
Let $\phi_\tau = \varphi'(F^\tau)u_2z_v$. Then the estimate \eqref{eqn:bilin} and Sobolev embedding imply
\begin{align*}
    \|\phi_\tau\|_{L^2(\mathcal{O})}  \leq & \|\varphi'(F^\tau)\|_{L^\infty(\mathcal{O})}\|u_2z_v\|_{L^2(\mathcal{O})} \leq c\|v\|_{L^2(\mathcal{O})},\\
    \|\nabla \phi_\tau\|_{L^2(\mathcal{O})} \leq & \|\nabla \varphi'(F^\tau )\|_{L^2(\mathcal{O})}\|u_2z_v\|_{L^\infty(\mathcal{O})}+  \|\varphi'(F^\tau)\|_{L^\infty(\mathcal{O})}\|\nabla(u_2z_v)\|_{L^2(\mathcal{O})}\\
    \leq & c(\|F_1\|_{H^1(\mathcal{O})}+\|F_2\|_{H^1(\mathcal{O})})\|v\|_{L^2(\mathcal{O})}.
\end{align*}
In sum, we obtain $\phi_\tau \in H_0^1(\mathcal{O})$, with 
$$
\|\phi_\tau \|_{H^1(\mathcal{O})}\leq c_2(\|F_1\|_{H^1(\mathcal{O})}+\|F_2\|_{H^1(\mathcal{O})})\|v\|_{L^2(\mathcal{O})}.
$$
By the Cauchy-Schwarz inequality, we have
\begin{align*}
    |(w,v)|&\leq \|F_2-F_1\|_{H^{-1}(\mathcal{O})}\int_0^1 \|\phi_\tau \|_{H^1(\mathcal{O})}\d \tau \\ 
    &\leq c_2(\|F_1\|_{H^1(\mathcal{O})}+\|F_2\|_{H^1(\mathcal{O})})\|F_2-F_1\|_{H^{-1}(\mathcal{O})}\|v\|_{L^2(\mathcal{O})}.
\end{align*}
Then taking supremum over $v\in L^2(\mathcal{O})$ with $\|v\|_{L^2(\mathcal{O})}\leq 1$ completes the proof of the lemma.
\end{proof}

Hence the forward map $\mathcal{G}$ satisfies Condition \ref{condreg} with  $k=1$. 
The stability result in \cite[Theorem 2.1 and Proposition 2.1]{Jin2023potential}, the uniform $H^2(\mathcal{O})$ bound on $\mathcal{G}(F)$ and Gargliardo-Nirenberg interpolation inequality \cite{BrezisMironescu:2018} yield that for all $M>0$ and $F\in H^1(\mathcal{O})$ with $\|F\|_{H^1(\mathcal{O})}\le M$,
$$
\|F-F_0\|_{L^2(\mathcal{O})} \lesssim \|\mathcal{G}(F)-\mathcal{G}(F_0)\|_{H^1(\mathcal{O})}^{\frac{1}{10}}\lesssim \|\mathcal{G}(F)-\mathcal{G}(F_0)\|_{L^2(\mathcal{O})}^{\frac{1}{20}}.
$$ 
This shows Condition \ref{condstab} with a H\"{o}lder exponent $\eta=\frac{1}{20}$.
For the Bayesian treatment, we again place a Gaussian prior $\Pi_{N,M}$ on $\Theta_M$, and use a Whittle-Mat\'{e}rn field with an RKHS $\mathcal{H}$ contained in $H^3(\mathcal{O})$. Fix $k=1$, $\alpha=3$, $\beta=1$ and $l=1$,
and set $b=-0.15,h=0.05,\rho=0.3$.  For each $d\in \{1, 2, 3\}$, one can choose some $\zeta>0$ (e.g., $\zeta=0.1$ for $d=1$, $\zeta=0.16$ for $d=2$, and $\zeta=0.25$ for $d=3$)
so that conditions \eqref{eqn:cond-P1}--\eqref{eqn:cond-P5} are satisfied.
Then Theorem~\ref{Thm:re-scaled consis} implies that the posterior  distribution $\Pi_{N,M}(\cdot\mid D_N)$ contracts around  $\mathcal{G}(F_0)$ at a rate $N^b = N^{-0.15}$.
Further, since the condition $\alpha + k \geq \frac{d(l+1)}{2}$ also holds, the contraction rate in Theorem~\ref{thm:VB-main} holds.

\begin{remark} 
Several existing results address the elliptic inverse potential problem \cite{doi:10.1137/18M1236137,nickl2018bernsteinvonmises}. In \cite{doi:10.1137/18M1236137}, the authors proved the existence and consistency of the MAP estimator under the assumptions $f\in H^\alpha(\mathcal O)$ with $\alpha>\frac{d}{2}$ and $g\in C^{\infty}(\mathcal O)$. In \cite{nickl2018bernsteinvonmises}, the Bernstein–von Mises theorem is established, where the true parameter satisfies $f\in C_c^s(\mathcal O)$ with $s>2+\frac{d}{2}$ and the boundary data $g$ satisfies $g\in C^{s+2}(\overline{\mathcal O})$.
In contrast, we treat the inverse potential problem in a low-regularity regime: the unknown coefficient  lies in $H^1(\mathcal O)$ and the source term $g$ lies in $L^2(\mathcal O)$. Using suitable continuity and stability arguments in weaker Sobolev norms, we derive the posterior contraction rates accordingly.
\end{remark}

\subsection{Inverse potential problem in subdiffusion}\label{ssec:subdiffusion}
Let $\mathcal{O}\subset \mathbb{R}^d$ ($d=1,2, 3$) be a smooth domain with a boundary $\partial\mathcal{O}$, and fix $T>0$ and $\beta=1$. Consider the following subdiffusion problem
\begin{equation}\label{eq:tf-forward}
\left\{
\begin{aligned}
\partial_t^{\alpha} u-\Delta u+fu &= h,
& & \text{in }\mathcal{O}\times(0,T),\\
u &= g,
& & \text{on }\partial\mathcal{O}\times(0,T),\\
u(\cdot,0) &= u_0,
& & \text{in }\mathcal{O},
\end{aligned}
\right.
\end{equation}
where $0<\alpha<1$ and $\partial_t^\alpha u(t)
:=\frac{1}{\Gamma(1-\alpha)}
\int_0^t (t-s)^{-\alpha} u'(s)\,\d s$ denotes the Caputo fractional
derivative \cite[p. 42]{jin2021fractional}.
The model is commonly employed to describe anomalously slow diffusion processes. Throughout, we make the following assumption on the problem data.

\begin{assumption}\label{ass:subdiff}
The following two conditions hold.
\begin{itemize}
\item[{\rm(i)}]  $0<m_0\leq u_0\in H^2(\mathcal O)$, $
0\leq h\in W^{1,4}(\mathcal O)$ and $
0<m_g\leq g\in H^{2}(\partial\mathcal O)$,
and $u_0=g$ on $\partial\mathcal{O}$;
\item[{\rm (ii)}] 
$ h + \Delta u_0 - \Lambda u_0 \ge 0$ in $\mathcal{O}$.
\end{itemize}
\end{assumption}

The compatibility condition $u_0=g$ in (i) on the boundary $\partial\mathcal{O}$ ensures good Sobolev regularity of the solution $u$, whereas (ii) ensures the nonnegativity of the Caputo derivative at the initial time.
We aim to recover the potential $f$ from the terminal data $\mathcal{G}(F):= u_f(\cdot,T)$ in the domain $\mathcal{O}$.
By the assumptions on  $u_0$ and $g$, we have \cite[Corollary 6.2]{jin2021fractional}:
\begin{equation}\label{eqn:subdif-apriori} 
\|u_f(\cdot,t)\|_{C([0,T];C(\overline{\mathcal{O}}))}\le U\quad\mbox{and}\quad
\|u_f(\cdot,t)\|_{C([0,T];H^1(\mathcal{O}))}\le c_1,
\end{equation}
where the constants $U$ and $c_1$ depend on $\lambda$, $\Lambda$, $\mathcal{O}$, $T$, $\alpha$, $u_0$, $h$ and $g$, but are independent of $f\in\mathcal{F}$. Thus there exists  $U < \infty$ such that $\sup_{F \in \Theta}\,\|\mathcal{G}(F)\|_{C(\overline{\mathcal{O}})} \le U.$ 

For any $f\in \mathcal{F}$, let $A_f$ denote the realization of the operator $-\Delta+f$ in $L^2(\mathcal{O})$ with a zero Dirichlet boundary condition. The operator $A_f$ satisfies the following property.  
\begin{lemma}\label{lem:ellipt-reg}
Let $f\in \mathcal{F}$ and $v\in H_0^1(\mathcal{O})$ satisfy $A_fv=g$. Then the following estimates hold
\begin{equation*}
    \|v\|_{H^2(\mathcal{O})}\leq c\|g\|_{L^2(\mathcal{O})}\quad \mbox{and} \quad \|v\|_{H^3(\mathcal{O})}\leq c(1+\|f\|_{H^1(\mathcal{O})})\|g\|_{H^1(\mathcal{O})},
\end{equation*}
where the constant $c$ depends only on $\lambda$, $\Lambda$ and $\mathcal{O}$.
\end{lemma}
\begin{proof}
The first estimate follows from the standard elliptic regularity theory. Next, by writing $A_fv=g$ as $A_0v=g-fv$. The first estimate and the Sobolev embedding $H^2(\mathcal{O})\hookrightarrow L^\infty(\mathcal{O})$ imply
\begin{align*}
    \|\nabla (g-fv)\|_{L^2(\mathcal{O})}&\leq \|\nabla g\|_{L^2(\mathcal{O})} +\|f\|_{L^\infty(\mathcal{O})}\|\nabla v\|_{L^2(\mathcal{O})}+\|\nabla f\|_{L^2(\mathcal{O})}\|v\|_{L^\infty(\mathcal{O})}\\
    & \leq c(\Lambda+\|f\|_{H^1(\mathcal{O})})\|g\|_{H^1(\mathcal{O})}.
\end{align*}
Then the elliptic regularity theory for the Dirichlet Laplacian $A_0$ implies the desired assertion.
\end{proof}

Then for any $\tilde g\in L^1(0,T;L^2(\mathcal O))$, the solution $w$ of
$$\partial_t^\alpha w+A_fw=\tilde g,\quad
\mbox{with } w(\cdot,0)=0,$$
can be represented as
\begin{align}\label{eq:solution_rep}
w(t)= \int_0^t E_f(t-s)\tilde g(s)\,\d s, 
\end{align}
with the solution operator
$E_f(t)$ given by 
$$
E_f(t):=t^{\alpha-1}E_{\alpha,\alpha}(-t^\alpha A_f),\quad t>0,$$
where $E_{\alpha,\alpha}(z)=\sum_{k=0}^\infty\frac{z^k}{\Gamma((k+1)\alpha)}$ denotes the Mittag--Leffler function \cite[Section 3.1]{jin2021fractional}.
The next lemma gives key smoothing properties of the operator $E_f(t)$ \cite[Theorem 6.4]{jin2021fractional}, where $\|\cdot\|_{\mathcal{L}(L^2(\mathcal{O})}$ denotes the operator norm from $L^2(\mathcal{O})$ into itself.
\begin{lemma}\label{lem:solution_operators}
For any $f\in\mathcal{F}$, $E_f(t)$ satisfies for $t\in(0,T]$:
\begin{align}
  \|A_f^\gamma E_f(t)\|_{\mathcal{L}(L^2(\mathcal{O}))}\le c t^{\alpha(1-\gamma)-1}, \quad \forall \gamma \in [0,1].
\end{align}
\end{lemma}

Furthermore, for any $f\in \mathcal{F}$, by the maximum principle, $u_{f}(T) \geq m_T>0$, with $m_T$ depending only on $m_0$, $m_g$, $T$, $\lambda$ and $\Lambda$. 
This result is similar to \cite[Lemma 2.3(ii)]{doi:10.1137/21M1446708}, except that the condition on $f$ is relaxed from $f\in \mathcal{F}\cap C(\overline{\mathcal{O}})$ to $f\in \mathcal{F}$. The proof of the lemma relies on a version of the maximum principle and is given in Appendix \ref{app:subdiff}. Note that the argument is substantially different from that in \cite{doi:10.1137/21M1446708}, and  requires a careful use of the weak formulation in \cite[Lemma 3.1]{Zacher:2008} and the convexity inequality in \cite[Corollary 4.1]{Zacher:2019}. 
\begin{lemma}\label{lem:max-princ}
Let Assumption \ref{ass:subdiff} hold, and let $u_f$ be the solution to problem \eqref{eq:tf-forward}. Then there exists  $m_T>0$, depending only on $m_0$, $m_g$, $T$, and $\Lambda$, such that $u_f(\cdot,T)\ge m_T>0$ in $\mathcal{O}$.
\end{lemma}

The solution map $\mathcal{G}$ is Lipschitz continuous in the following sense.
\begin{lemma}\label{ex3: Lip}
There exists $L>0$ such that
\begin{equation*}\label{eq:ex3Lip}
\bigl\|\mathcal{G}(F_1)-\mathcal{G}(F_2)\bigr\|_{L^{2}(\mathcal{O})}
\leq  
L(\|F_1\|_{H^1(\mathcal{O})}+\|F_2\|_{H^1(\mathcal{O})})\|F_2-F_1\|_{H^{-1}(\mathcal{O})},\quad \forall F_1,F_2\in \Theta.
\end{equation*}
\end{lemma}
\begin{proof}
Let $u_i = \mathcal{G}(F_i)$, $f_i=\varphi\circ F_i$, $i=1,2$ and let $w:=u_1-u_2$. Then $w$ satisfies $$\partial_t^\alpha w+A_{f_1}  w=(f_2-f_1)u_2,
\quad \mbox{with } w(\cdot,0)=0.$$
By the representation \eqref{eq:solution_rep}, for any $v\in L^2(\mathcal{O})$, there holds 
$$
(w(T),v)=\bigg(\int_0^T E_{f_1}(T-s)\bigl[(f_2-f_1)u_2(\cdot,s)\bigr]\,\d s,v\bigg) =\int_0^T \bigl(f_2-f_1,\,u_2(\cdot,s)\zeta_v(s)\bigr)\,\d s,
$$
with $\zeta_v(s):=E_{f_1}(T-s)v$. 
For the term $f_1-f_2$, we employ the identity \eqref{eqn:F-fdc}. Let
$a_\tau(x):=\varphi'\bigl(F_1+\tau(F_2-F_1)\bigr)(x)$ and $\phi_{\tau,s,v}(x):=a_\tau(x)\,u_2(x,s)\,\zeta_v(x,s)$. 
Note that  $$\|a_\tau\|_{L^\infty(\mathcal{O})}\le \|\varphi'\|_{L^\infty(\mathbb{R})}\quad \mbox{and} 
\quad \|\nabla a_\tau\|_{L^2(\mathcal{O})}
\le \|\varphi''\|_{L^\infty(\mathbb{R})}\bigl(\|F_1\|_{H^1(\mathcal{O})}+\|F_2\|_{H^1(\mathcal{O})}\bigr).$$
Fix $\varepsilon\in(0,1)$ such that $1+\varepsilon>\frac{d}{2}$. By Lemma~\ref{lem:solution_operators}, for any $s\in(0,T)$, there holds
\begin{align}\label{eqn:bdd-zeta}  \|\zeta_v(\cdot,s)\|_{H^{1+\varepsilon}(\mathcal{O})}
\le\|A_{f_1}^{\frac{1+\varepsilon}{2}} E_{f_1}(T-s) v\|_{L^2(\mathcal{O})} \le c(T-s)^{\frac{1-\varepsilon}{2}\alpha-1}\ \|v\|_{L^2(\mathcal{O})}.
\end{align}
The continuous Sobolev embedding $H^{1+\varepsilon}(\mathcal{O})\hookrightarrow L^\infty(\mathcal{O})$, the \textit{a priori} estimates in \eqref{eqn:subdif-apriori}, and the identity $\nabla\phi_{\tau,s,v}=(\nabla a_\tau)\,u_2\zeta_v+a_\tau\nabla(u_2\zeta_v)$ imply
\begin{align*}
\|\phi_{\tau,s,v}\|_{L^2(\mathcal{O})}
&\le \|a_\tau\|_{L^\infty(\mathcal{O})}\|u_2\|_{L^\infty(\mathcal{O})}\|\zeta_v\|_{L^2(\mathcal{O})}
\le cU\|\varphi'\|_{L^\infty(\mathbb{R})} \|\zeta_v\|_{H^{1+\varepsilon}(\mathcal{O})},\\  
\|\nabla\phi_{\tau,s,v}\|_{L^2(\mathcal{O})}
&\le \|\nabla a_\tau\|_{L^2(\mathcal{O})}\,\|u_2\|_{L^\infty(\mathcal{O})}\|\zeta_v\|_{L^\infty(\mathcal{O})}
+\|a_\tau\|_{L^\infty(\mathcal{O})}\,\|\nabla(u_2\zeta_v)\|_{L^2(\mathcal{O})}\\
&\le  \|\nabla a_\tau\|_{L^2(\mathcal{O})}( Uc\|\zeta_v\|_{H^{1+\varepsilon}(\mathcal{O})})+\|a_\tau\|_{L^\infty(\mathcal{O})}\,(c_1\,c\|\zeta_v\|_{H^{1+\varepsilon}(\mathcal{O})} + U\,\|\zeta_v\|_{H^1(\mathcal{O})}).
\end{align*}
Hence for some $c=c(\mathcal{O}, d, c_1, U,\lambda)$, there holds
\begin{equation*}
\|\phi_{\tau,s,v}\|_{H^1(\mathcal{O})}
\le c\bigl(\|F_1\|_{H^1(\mathcal{O})}+\|F_2\|_{H^1(\mathcal{O})}\bigr)
\|\zeta_v(\cdot,s)\|_{H^{1+\varepsilon}(\mathcal{O})}.
\end{equation*}
Thus $\phi_{\tau,s,v}\in H_0^1(\mathcal{O})$.
By the Cauchy-Schwarz inequality and the estimate \eqref{eqn:bdd-zeta}, we deduce
\begin{align*}
&|(w(\cdot,T),v)|
\le \|F_2-F_1\|_{H^{-1}(\mathcal{O})}\int_0^T\int_0^1 \|\phi_{\tau,s,v}\|_{H^1(\mathcal{O})}\ \d \tau\,\d s\\
\le &c\bigl(\|F_1\|_{H^1(\mathcal{O})}+\|F_2\|_{H^1(\mathcal{O})}\bigr)
\|F_2-F_1\|_{H^{-1}(\mathcal{O})}\|v\|_{L^2(\mathcal{O})}
\int_0^T (T-s)^{\alpha-1-\frac{\alpha(1+\varepsilon)}{2}}\d s\\
\le &L\bigl(\|F_1\|_{H^1(\mathcal{O})}+\|F_2\|_{H^1(\mathcal{O})}\bigr)
\|F_2-F_1\|_{H^{-1}(\mathcal{O})}\|v\|_{L^2(\mathcal{O})}.
\end{align*}
Taking the supremum over all $v\in L^2(\mathcal{O})$ with $\|v\|_{L^2(\mathcal{O})}\leq 1$ gives the desired estimate. 
\end{proof}

Hence the forward map $\mathcal{G}$ satisfies Condition \ref{condreg} with $k=1$. The next result gives a conditional stability estimate, thereby showing Condition \ref{condstab}. \begin{lemma}\label{ex3:cond_stab}
Let Assumption \ref{ass:subdiff} hold. For any $ F_1,F_2\in \Theta_M$, with $f_i:=\varphi\circ F_i\in\mathcal{F}$, $i=1,2$, there exist $T_0>0$, $c_{\rm s}=c_{\rm s}(M, T)$ and some $\eta\in (0,1)$ such that for all $T \geq T_0$,

$$
\left\|F_1-F_2\right\|_{L^2(\mathcal{O})} \leq c_{\rm s}\left\|\mathcal{G}(F_1)-\mathcal{G}\left(F_2\right)\right\|_{L^2(\mathcal{O})}^\eta.
$$    
\end{lemma}
\begin{proof}
Let $w=u_1-u_2$. At $t=T$, we have $(f_2- f_1)u_2(T) = -\partial_t^\alpha w(T) - \mathcal{A}_{f_1} w(T)$, and by Lemma \ref{lem:max-princ}, we deduce 
\begin{align}\label{eq:basic_stab}
\|f_2-f_1\|_{L^2(\mathcal{O})}
\le
m_T^{-1}\bigl(\|\partial_t^\alpha w(T)\|_{L^2(\mathcal{O})}
+\|w(T)\|_{H^2(\mathcal{O})}\bigr).
\end{align}
Further, by \cite[Lemma 3.2]{doi:10.1137/21M1446708}, for any $\epsilon\in\left(0,\min(1,2-\frac d2)\right)$, there exists $T_0>0$ such that for all $T\ge T_0$, there holds
\begin{align}\label{eq:dtalpha_decay}
\|\partial_t^\alpha w(T)\|_{L^2(\mathcal{O})}
\le c \max \bigl(T^{-\alpha},\, T^{-(1-\epsilon)\alpha}\bigr)
\|f_2-f_1\|_{L^2(\mathcal{O})}.
\end{align}
The estimates \eqref{eq:dtalpha_decay} and \eqref{eq:basic_stab} together yield that for large $T$,
\begin{align}\label{eq:ex3_F_by_H2}
\|F_1-F_2\|_{L^2(\mathcal O)}\le c_{\rm s}\,\|w(T)\|_{H^2(\mathcal O)},
\end{align}
with $c_{\rm s} = c_{\rm s}(M,T)$.
Let $\tilde g(s)=(f_2-f_1)u_2(s)$. Following the argument of Lemma \ref{lem:ellipt-reg}, the estimates in \eqref{eqn:subdif-apriori} imply 
\begin{equation}\label{eqn:bdd-tildeg}
    \|\tilde g\|_{L^\infty(0,T;H^1(\mathcal O))}\leq c(\|f_1-f_2\|_{H^1(\mathcal O)}+\|f_1-f_2\|_{L^\infty(\mathcal{O})})\leq c, 
\end{equation}
where $c$ depends on $U$, $c_1$, $\lambda$ and $\Lambda$. 
For any $\rho\in(0,1)$, by Lemma \ref{lem:ellipt-reg}, the $H^{2+\rho}(\mathcal{O})$ norm of $w(T)$ is bounded by $\|A_{f_1}^{1+\frac\rho2} w(T)\|_{L^2(\mathcal{O})}$. By the representation \eqref{eq:solution_rep}, Lemma \ref{lem:solution_operators} with $\gamma =\frac{1}{2} + \frac{\rho}{2}$ and the estimate \eqref{eqn:bdd-tildeg}, we derive
\begin{align*}
\|w(T)\|_{H^{2+\rho}(\mathcal{O})} &\leq \int_0^T \left\| A_{f_1}^{\frac12+\frac\rho2} E_{f_1}(-\tau) A_{f_1}^\frac12\tilde g(T-\tau) \right\|_{L^2(\mathcal{O})} \, \d \tau\\   
 &\lesssim \int_0^T \tau^{(\frac12-\frac\rho2)\alpha-1} \| A_{f_1}^\frac12 \tilde g(T-\tau) \|_{L^2(\mathcal{O})} \, \d \tau\\
&=  \|\tilde g\|_{L^\infty(0,T; H^1(\mathcal{O}))}\int_0^T \tau^{-1 + \alpha(\frac12 - \frac\rho2)} \, \d \tau<\infty.
\end{align*}
Using the estimate \eqref{eq:ex3_F_by_H2} and the Gagliardo–Nirenberg interpolation inequality \cite{BrezisMironescu:2018}
$$\|w(T)\|_{H^2(\mathcal{O})}
\le c,\|w(T)\|_{L^2(\mathcal{O})}^{\eta}\,\|w(T)\|_{H^{2+\rho}(\mathcal{O})}^{1-\eta},$$ with $\eta=\frac{\rho}{2+\rho}$, we obtain the desired assertion. 
\end{proof}

The Bayesian analysis follows identically to that in Section~\ref{ssec:potential}.  Using the Whittle--Mat\'{e}rn prior and the parameter choice therein, 
Theorems~\ref{Thm:re-scaled consis} and~\ref{thm:VB-main} are applicable and yield a posterior contraction rate  $N^{-0.15}$. Since the solution map satisfies the stability assumption with $\eta=\frac{\rho}{2+\rho}$, choosing $\rho$ sufficiently close to 1 yields a stability exponent $\eta$ close to $\frac{1}{3}$, and therefore a posterior contraction rate for the unknown parameter is close to $N^{-0.05}$.

\begin{remark}
Kow and Wang \cite{kow2025consistency} investigated the recovery of the potential in a subdiffusion model involving the Riemann-Liouville fractional derivative from noisy observations of the solution at randomly sampled space--time points $(t_i,x_i)\in \mathcal{O}\times(0,T)$. Note that the Riemann-Liouville fractional derivative is the closure of the Caputo one in a suitable Hilbert space. In contrast, we investigate the subdiffusion model with the Caputo fractional derivative and focus on the recovery of the potential from the terminal observation over the spatial domain $\mathcal{O}$. These different modeling choices lead to different forward operators and require different analytical treatments in the posterior contraction rate analysis. The present analysis allows for a lower regularity setting than that in \cite{kow2025consistency}.
\end{remark}

\section*{Acknowledgments} The authors are grateful to two anonymous referees for their constructive comments which have led to a significant improvement in the quality of the paper.

\appendix

\section{Bound on the contraction rate}

In the proof of Theorem \ref{Thm:re-scaled consis}, we have introduced the system of constraints \eqref{eqn:cond-P1}--\eqref{eqn:cond-P5} for  determining the  contraction rate $\delta_N = N^b$ ($b<0$). For the problem setting $d\geq1$, $l\geq 0$, $k\geq0$, and $\alpha> \beta+\frac{d}{2}$ in Section \ref{ssec:Proof Bayes case}, we derive a sharp lower bound of $b$ ensuring the existence of a feasible triple $(\zeta, \rho, h)$ with $\zeta>0$, $\rho\geq 0$, and $h>0$.
\begin{proposition}\label{prop:rate}
Let
$$
b_1:= \frac{k+\beta}{2(k+\alpha)}\quad\mbox{and} \quad b_2:= \frac{\alpha+k}{d+2(\alpha+k) + \frac{(1+l)d(\alpha-\beta)}{k+\beta}}.
$$
The system of constraints \eqref{eqn:cond-P1}--\eqref{eqn:cond-P5} admits a feasible solution $(\zeta, \rho, h)$ with $\zeta>0$, $\rho\geq 0$, and $h>0$ if and only if 
$$-\min\{b_1, b_2\}  <b < 0.$$
A sharp bound on the contraction rate in Theorem \ref{Thm:re-scaled consis} is given by $b_* = -\min\{b_1, b_2\}$.
\end{proposition}
\begin{proof}
First we reduce the system by eliminating $h$. From \eqref{eqn:cond-P4}, we have $h = b + \frac{1}{2} - \rho$. Substituting this into \eqref{eqn:cond-P2}, \eqref{eqn:cond-P3} and \eqref{eqn:cond-P5} yields the following equivalent conditions on $\zeta$ and $\rho$:
\begin{align}
\rho &\geq \zeta \frac{\alpha-\beta}{d}, \label{eq:red_P2} \\
\rho &\leq  b + \frac{(1+2b)(\alpha+k)}{d}, \label{eq:red_P3} \\
\rho &\leq \frac{1}{1+l}\left(b + \frac{(1+2b)(\alpha+k)}{d}\right). \label{eq:red_P5}
\end{align}
Let $C(b):=b + \frac{(1+2b)(\alpha+k)}{d}$. Since $l \geq 0$, condition \eqref{eq:red_P5} is strictly sharper than or equal to \eqref{eq:red_P3}. In view of \eqref{eq:red_P2}, \eqref{eq:red_P5}, and the positivity requirement $h > 0$, an admissible $\rho$ exists if and only if the intersection interval
$$I_\rho := \left[ \zeta\frac{\alpha-\beta}{d}, \, \min\left\{ b+\frac12, \, \frac{C(b)}{1+l} \right\} \right)$$
is nonempty.
Now, the constraint \eqref{eqn:cond-P1} requires $\zeta > -\frac{bd}{k+\beta}$, which together with $\alpha > \beta$  gives  a strict lower bound on the left endpoint of $\mathcal{I}_\rho$:
$$
\zeta\frac{\alpha-\beta}{d} > -b\frac{\alpha-\beta}{k+\beta}.
$$
Therefore, a feasible pair $(\zeta, \rho)$ exists if and only if $-b\frac{\alpha-\beta}{k+\beta}$ is strictly bounded from above by the right endpoint of $I_\rho$, which leads to the following coupled inequalities:
\begin{align*}
-b\frac{\alpha-\beta}{k+\beta} &< b+\frac12 \iff b > -b_1,\\
-b\frac{\alpha-\beta}{k+\beta} &< \frac{C(b)}{1+l} \iff b > -b_2.
\end{align*}
Combining these two requirements yields the desired assertion $b > -\min\{b_1, b_2\}$. 

\noindent Now we prove the sufficiency. If $-\min\{b_1, b_2\} < b < 0$, then there exists $\epsilon > 0$ sufficiently small and $\zeta = -\frac{bd}{k+\beta} + \epsilon$ such that $\zeta\frac{\alpha-\beta}{d} < \min\left\{ b+\frac12, \frac{C(b)}{1+l} \right\}$. Any choice of $\rho$ inside the  nonempty interval $I_\rho$ alongside $h = b + \frac{1}{2} - \rho > 0$ satisfies conditions  \eqref{eqn:cond-P1}--\eqref{eqn:cond-P5} simultaneously.
\end{proof}

\begin{remark}
The two thresholds $b_1$ and $b_2$ in Proposition \ref{prop:rate} reflect two distinct mechanisms: prior oversmoothing and PDE stability. The first threshold $b_1$ captures the oversmoothing property of the prior and is consistent with Castillo's Gaussian-prior contraction lower bound \cite[Theorem 2]{Castillo:2008}, after translating the $(H^{k}(\mathcal{O}))^*$-loss into an $L^2(\mathcal{O})$-loss at the spectral-sequence level.
The second threshold, $b_2$, arises from the nonlinearity of the forward operator. Since larger values of $l$ correspond to weaker stability properties of the PDE model, the admissible contraction exponent becomes smaller. This captures the expected degradation of the posterior contraction rate when the forward operator has poorer stability.
\end{remark}

\section{Auxiliary estimates}

First we give an elliptic regularity estimate.
\begin{lemma}\label{lem:difference_quotient}
Let $\mathcal O\subset\mathbb R^2$ be a smooth bounded domain, and let $a\in H^{1+\varepsilon}(\mathcal O)$ satisfy
$$
0<\lambda\le a(x)\le \Lambda \quad \text{a.e. in }\mathcal O\quad \text{and}\quad \|a\|_{H^{1+\varepsilon}(\mathcal O)} \le M .
$$
Then for every $0<s<\varepsilon<1$,
there exists $c=c(\lambda,\Lambda,\mathcal O,s,M)$
such that for every $v\in L^2(\mathcal O)$, the weak solution
$z\in H_0^1(\mathcal O)$ of
$$
-\nabla\cdot(a\nabla z)=v, \quad \text{in }\mathcal O
$$
satisfies
\begin{equation}\label{eqn:ell-reg-frac}
\|\nabla z\|_{H^s(\mathcal O)}\le c\|v\|_{L^2(\mathcal O)}.
\end{equation}
\end{lemma}
\begin{proof}
The proof of the regularity estimate \eqref{eqn:ell-reg-frac} follows the difference quotient approach \cite[Section 6.3]{evans2010partial} and is lengthy and divided into three steps. We use extensively the $L^2(\mathcal{O})$ inner product $(\cdot,\cdot)$, and the $L^2(\mathcal{K})$ inner product $(\cdot,\cdot)_{L^2(\mathcal{K})}$ .

\medskip
\noindent 
\underline{Step 1}. \label{lem:step1}
Since $\mathcal O\subset\mathbb R^2$ and $a\in H^{1+\varepsilon}(\mathcal O)$, the Sobolev embedding theorem yields for $0<\gamma<\varepsilon$,
$$a\in C^{0,\gamma}(\overline{\mathcal O})
\quad\mbox{and}\quad
\|a\|_{C^{0,\gamma}(\overline{\mathcal O})}
\le c \|a\|_{H^{1+\varepsilon}(\mathcal O)}.$$
For any $\mathcal{U}\Subset \mathcal O$, we choose a cut-off function $\eta\in C_c^\infty(\mathcal O)$ such that $0\le \eta\le 1$ on $\mathcal O$ and $\eta\equiv 1$ on $\mathcal{U}$, and denote $\operatorname{supp}(\eta)$ by $\mathcal{K}$. For any function $w$ defined on $\mathcal O$, and for $|h|$ sufficiently small, we define the finite difference in the
$e_k$-direction, $k=1,2$, by
$$
\delta_h^k w(x):=w(x+he_k)-w(x), \quad\forall x\in \mathcal O_h, 
$$
with $\mathcal O_h:=\{x\in\mathcal O:\ x+he_k\in\mathcal O\}.$ The following argument does not depend on any special property of $e_k$, and the estimate holds for any unit direction $e\in\mathbb{S}^1$. Let $\delta_h=\delta_h^k$, and by taking $|h|<h_0:=\frac{1}{2}\operatorname{dist}(\mathcal{K},\partial\mathcal{O})$, we can use $\phi:=\delta_{-h}(\eta^2 \delta_h z)\in H_0^1(\mathcal O)$ as a test function in the weak formulation of $z$:
\begin{equation}\label{def:test}
(a\nabla z,\nabla\delta_{-h}(\eta^2\delta_h z))
=(v,\delta_{-h}(\eta^2\delta_h z)). 
\end{equation}
By the discrete integration by parts formula \cite[Section 5.8.2]{evans2010partial}, we obtain
$$(a\nabla z,
\nabla\delta_{-h}(\eta^2\delta_h z))
=(\delta_h(a\nabla z), \nabla(\eta^2\delta_h z)).$$
For any $x\in\mathcal O_h$, the following product rule for the finite difference holds
$$
\delta_h(a\nabla z)= (\tau_h a)\nabla(\delta_h z)+(\delta_h a)\nabla z,
\quad \text{with }\tau_h a(x):=a(x+he_k),
$$
and by the product rule, there also holds $\nabla(\eta^2\delta_h z)
=
\eta^2\nabla(\delta_h z)
+
2\eta\,\delta_h z\,\nabla\eta$. Hence, we obtain
\begin{align*}
(\eta^2 (\tau_h a),|\nabla(\delta_h z)|^2)
=&
-(2\eta (\tau_h a)\,\delta_h z\, \nabla(\delta_h z), \nabla\eta)
-(\eta^2 (\delta_h a)\,\nabla z, \nabla(\delta_h z)) \\
&
-(2\eta (\delta_h a)\,\delta_h z\, \nabla z, \nabla\eta)
+(v,\delta_{-h}(\eta^2\delta_h z)):=\sum_{i=1}^4{\rm I}_i.
\end{align*}
Now we estimate the four terms $({\rm I}_{i})_{i=1}^4$ separately. Using the estimate \cite[pp. 277-278, Section 5.8.2]{evans2010partial} 
$$\|\delta_h z\|_{L^2(\mathcal{K})}\le |h|\,\|\nabla z\|_{L^2(\mathcal{K})},$$ 
the estimate $\|\tau_h a\|_{L^\infty(\mathcal{K})}
\le \|a\|_{L^\infty(\mathcal O)}
\le \Lambda$ and Young's inequality, we have
\begin{align*}
\left|{\rm I}_1\right|
&\le c(|\delta_h z|,\eta |\nabla(\delta_h z)|)_{L^2(\mathcal{K})} \le c \|\delta_h z\|_{L^2(\mathcal{K})}\,\|\eta \nabla(\delta_h z)\|_{L^2(\mathcal{K})} \\
&\leq \tfrac{\lambda}{4}\|\eta\nabla(\delta_h z)\|_{L^2(\mathcal{O})}^2 + c|h|^2\|\nabla z\|_{L^2(\mathcal O)}^2.
\end{align*}
Since \(a\in C^{0,\gamma}(\overline{\mathcal O})\),
we have $\|\delta_h a\|_{L^\infty(\mathcal{K})}\le c|h|^\gamma\|a\|_{C^{0,\gamma}(\mathcal O)}$. Hence 
\begin{align*}
|{\rm I}_2|
&\le \|\delta_h a\|_{L^\infty(\mathcal{K})} (\eta^2|\nabla z|,|\nabla(\delta_h z)|)
\le c|h|^\gamma \|\nabla z\|_{L^2(\mathcal O)}\,\|\eta\nabla(\delta_h z)\|_{L^2(\mathcal O)}.
\end{align*}
Following the argument of \cite[Section 5.8.2]{evans2010partial}, we can similarly derive 
\begin{align*}
\left|{\rm I}_3\right|
&\le c\|\delta_h a\|_{L^\infty(\mathcal{K})} \|\delta_h z\|_{L^2(\mathcal{K})}\,\|\nabla z\|_{L^2(\mathcal{K})} \le c|h|^{1+\gamma}\|\nabla z\|_{L^2(\mathcal{K})}^2.
\end{align*}
Last, by \cite[Section 5.8.2]{evans2010partial} and the duality argument, we have the estimate $\|\delta_h v\|_{H^{-1}(\mathcal{K})}\le c|h|\|v\|_{L^2(\mathcal O)}$ and deduce
\begin{align*}
\left|{\rm I}_4\right|&=\left|(\delta_h v,\eta^2\delta_h z)\right|
\le c|h|\,\|v\|_{L^2(\mathcal O)}\|\nabla(\eta^2\delta_h z)\|_{L^2(\mathcal O)} \\
&\le c|h|\,\|v\|_{L^2(\mathcal O)}\big(\|\eta\nabla(\delta_h z)\|_{L^2(\mathcal O)}+\|\delta_h z\|_{L^2(\mathcal O)}\big) \\
&\le c|h|\,\|v\|_{L^2(\mathcal O)}\|\eta\nabla(\delta_h z)\|_{L^2(\mathcal O)} + c|h|^2\|v\|_{L^2(\mathcal O)}\|\nabla z\|_{L^2(\mathcal O)}.
\end{align*}
Substituting the bounds on $({\rm I}_i)_{i=1}^4$ into \eqref{def:test} and using the standard energy estimate $\|\nabla z\|_{L^2(\mathcal O)}\le c\|v\|_{L^2(\mathcal O)}$ \cite[Section 6.2.2]{evans2010partial}, we derive
$$
\lambda \|\eta\nabla(\delta_h z)\|_{L^2(\mathcal O)}^2
\le \tfrac{\lambda}{4}\|\eta\nabla(\delta_h z)\|_{L^2(\mathcal O)}^2
+ c|h|^\gamma \|v\|_{L^2(\mathcal O)}\|\eta\nabla(\delta_h z)\|_{L^2(\mathcal O)}
+ c|h|^{1+\gamma}\|v\|_{L^2(\mathcal O)}^2.
$$
Collecting the terms and then solving the quadratic inequality of $\|\eta\nabla(\delta_h z)\|_{L^2(\mathcal O)}$ yield
$$
\|\eta\nabla(\delta_h z)\|_{L^2(\mathcal O)} \le c|h|^\gamma \|v\|_{L^2(\mathcal O)}.
$$
Since $0\le \eta\le 1$ on $\mathcal O$, this estimate gives
\begin{equation}\label{eq:difference_estimate}
\|\delta_h(\nabla z)\|_{L^2(\mathcal{O}_h)} \le c|h|^\gamma \|v\|_{L^2(\mathcal O)}.
\end{equation}
By the difference quotient characterization of fractional Sobolev spaces \cite[Exercise 17.23]{leoni2017firstSobolev} and note that $\mathcal{U}\Subset\mathcal{O}_h$, we obtain for every $0<s<\gamma$,
\begin{align}
\|\nabla z\|_{H^s(\mathcal{U})}^2
&\lesssim \|\nabla z\|_{L^2(\mathcal O)}^2
+\int_0^{h_0}
\frac{\|\delta_h \nabla z\|_{L^2(\mathcal O_h)}^2}{h^{1+2s}}\,\d h\label{eqn:frac-Sobolev}\\
&\lesssim \|v\|_{L^2(\mathcal O)}^2+\int_0^{h_0}
\frac{h^{2\gamma}\|v\|_{L^2(\mathcal O)}^2}{h^{1+2s}}\,\d h \lesssim\|v\|_{L^2(\mathcal O)}^2.\nonumber
\end{align}
The argument remains valid for any symmetric and uniformly elliptic $A(x) \in C^{0,\gamma}(\overline{\mathcal{O}};\mathbb{R}^{2 \times 2})$: 
$$
\lambda|\xi|^2 \leq A(x) \xi \cdot \xi \leq \Lambda|\xi|^2, \quad \forall \xi \in\mathbb{R}^2, \ x\in\mathcal O.
$$

\noindent \underline{Step 2}.
Now we derive the boundary estimate. Fix $x_0\in \partial\mathcal O$. For $r>0$, set
$B_r:=\{y\in\mathbb{R}^2:|y|<r\}$ and
$B_r^+:=B_r\cap\{y_2>0\}$. Since $\partial\mathcal O$ is smooth, there exist a
neighborhood $\mathcal{W}$ of $x_0$ and a $C^{2}$-diffeomorphism 
$\Phi:B_2\to \mathcal{W}$ such that
$\Phi(B_2^+)=\mathcal{W}\cap\mathcal O$ and $\Phi(B_2\cap\{y_2=0\})=\mathcal{W}\cap\partial\mathcal O$.
Let
$$
\widehat A(y)
=|\det D\Phi(y)|\,a(\Phi(y))\,D\Phi(y)^{-1}D\Phi(y)^{-T},\quad \widehat z(y)= z(\Phi(y)),\quad \widehat v(y)=|\det D\Phi(y)| v(\Phi(y)).
$$
Then $\widehat z$ solves in $B_2^+$ the following equation
$$
-\nabla_y\cdot(\widehat A(y)\nabla_y \widehat z)=\widehat v,\quad\mbox{in } B_2,
$$
with a zero Dirichlet boundary condition $\widehat z=0$ on $\{y_2=0\}$, where
$\widehat A\in C^{0,\gamma}(\overline{B_2^+};\mathbb{R}^{2\times 2})$ is symmetric and uniformly elliptic,
with bounds depending only on $\lambda,\Lambda,M $ and the chart $\Phi$ \cite[Section 6.3.2]{evans2010partial}. Moreover, we take $\Phi$ to be a arc-length-normal chart on the flat boundary, so that $\widehat A_{12}(y_1,0)=\widehat A_{21}(y_1,0)=0$.
Now define the odd reflections of $\widehat z$ and $\widehat v$ across $\{y_2=0\}$ by 
$\widetilde z(y_1,-y_2):=-\widehat z(y_1,y_2)$ and $
\widetilde v(y_1,-y_2):=-\widehat v(y_1,y_2)$,
with $\widetilde z \in H^1(B_2)$, and $\widetilde v \in L^2(B_2)$, and extend the coefficient matrix $\widehat A(y)= \big(\widehat A_{ij}(y)\big)_{1\le i,j\le 2}$ by 
\begin{align*}
&\widetilde A_{11}(y_1,-y_2)=\widehat A_{11}(y_1,y_2),
\quad
\widetilde A_{22}(y_1,-y_2)=\widehat A_{22}(y_1,y_2),\\
&\widetilde A_{12}(y_1,-y_2)=-\widehat A_{12}(y_1,y_2),
\quad
\widetilde A_{21}(y_1,-y_2)=-\widehat A_{21}(y_1,y_2).
\end{align*}
Then the extension $\widetilde A\in C^{0,\gamma}(\overline{B_2};\mathbb{R}^{2\times 2})$ is symmetric and uniformly elliptic. Then the function $\widetilde z \in  H^1(B_2)$ is a weak solution of
$$
-\nabla\cdot(\widetilde A\nabla \widetilde z)=\widetilde v.
$$
By applying the estimate \eqref{eq:difference_estimate} to $\widetilde z$ in the ball $B_2$ and take $B_{\frac32}\Subset B_2$,
we obtain 
\begin{equation*}\label{eq:difference_estimate_boundary}
\|\delta_h(\nabla \widetilde z)\|_{L^2(B_{\frac32})} \le c|h|^\gamma \|v\|_{L^2(\mathcal O)}.
\end{equation*}
By the difference quotient characterization \eqref{eqn:frac-Sobolev} of fractional Sobolev spaces,
this estimate implies $$\|\nabla_y \widehat z\|_{H^s(B_1^+)}
\le c\|v\|_{L^2(\mathcal O)},\quad\mbox{for }
0<s<\gamma.$$
Since the $H^s$-norm is preserved up to equivalence under $C^2$-diffeomorphisms, this implies the same estimate for $\nabla z$ in a
neighborhood $\mathcal{W}_{x_0}:=\Phi(B_1^+)$ of $x_0$:
\begin{equation}\label{lem:boundary_neighbourhood}
   \|\nabla_x z\|_{H^s(\mathcal{W}_{x_0}\cap\mathcal O)}
\le c\|v\|_{L^2(\mathcal O)},
\end{equation}
with the constant depending only on
$\lambda$, $\Lambda$, $\mathcal O$, $M$ and the chart $\Phi$.

\noindent \underline{Step 3}.
Since the boundary $\partial\mathcal{O}$ is compact, we can extract a finite cover $\{\mathcal{W}_{x_i}\}_{i=1}^N$ from the collection $\{\mathcal{W}_x\}_{x \in \partial\mathcal{O}}$ such that $\partial\mathcal{O} \subset \bigcup_{i=1}^N \mathcal{W}_{x_i}$. Let $\mathcal{V} = \bigcup_{i=1}^N\mathcal{W}_{x_i}$. We combine the estimates in \eqref{lem:boundary_neighbourhood} and obtain $ \|\nabla z\|_{H^s(\mathcal{V} \cap\mathcal O)}
\le c\|v\|_{L^2(\mathcal O)}$. Since $\overline{\mathcal{O}} \setminus \mathcal{V} \Subset \mathcal O$, we can fix a compact set $\mathcal{U}$ with $\overline{\mathcal{O}} \setminus \mathcal{V} \Subset \mathcal{U} \Subset \mathcal O$. By the interior estimate from Step 1, we have $$\|\nabla z\|_{H^s(\mathcal{U})}
\le c\|v\|_{L^2(\mathcal O)},$$
with $c=c(\lambda,\Lambda,\mathcal O,M)$. 
By the characterization of Sobolev norms \cite[(2.2)]{di2012guidefractional}, $$\|\nabla z\|_{H^s(\mathcal O)}^2 \simeq \|\nabla z\|_{L^2(\mathcal O)}^2+[\nabla z]_{H^s(\mathcal O)}^2
,\quad\mbox{with }
[\nabla z]_{H^s(\mathcal O)}^2
=
\int_{\mathcal O}\int_{\mathcal O}
\frac{|\nabla z(x)-\nabla z(y)|^2}{|x-y|^{2+2s}}\d x\d y.
$$
Since $\overline{\mathcal O}\setminus \mathcal V \Subset \mathcal{U}$, we can choose $r>0$ small enough such that $\operatorname{dist}\bigl(\overline{\mathcal O}\setminus \mathcal V, \overline{\mathcal O}\setminus \mathcal{U}\bigr)>2r$.
Then, if $x\in \overline{\mathcal O}\setminus \mathcal V$ and
$|x-y|<r$, we have $y\in \mathcal{U}$. Similarly, if
$x\in \overline{\mathcal O}\setminus\mathcal{U}$ and $|x-y|<r$, then
$x,y\in \mathcal V$. Hence the near-field part of the fractional seminorm
is controlled by the local seminorms on $\mathcal{U}$ and $\mathcal V$:
$$\int_{\mathcal O}\int_{\mathcal O\cap B_r(x)}
\frac{|\nabla z(x)-\nabla z(y)|^2}{|x-y|^{2+2s}}\d y\d x
\le
c\Bigl(
[\nabla z]_{H^s(\mathcal{U})}^2+
[\nabla z]_{H^s(\mathcal V)}^2
\Bigr).$$
For the far-field part, since $|x-y|\ge r$, the integral is bounded by
$$\int_{\mathcal O}\int_{\mathcal O\setminus B_r(x)}
\frac{|\nabla z(x)-\nabla z(y)|^2}{|x-y|^{2+2s}}\d y\d x
\le
c_r\|\nabla z\|_{L^2(\mathcal O)}^2.$$
Consequently,
$$[\nabla z]_{H^s(\mathcal O)}^2
\le
c\Bigl(
[\nabla z]_{H^s(\mathcal{U})}^2+
[\nabla z]_{H^s(\mathcal V)}^2+
\|\nabla z\|_{L^2(\mathcal O)}^2\Bigr).$$
Hence \
$$\|\nabla z\|_{H^s(\mathcal O)}
\le c(\|\nabla z\|_{H^s(\mathcal{U})}+\|\nabla z\|_{H^s(\mathcal V)})\leq C\|v\|_{L^2(\mathcal O)}.$$
This completes the proof of the lemma.
\end{proof}

Next we give a bound on the product of a Sobolev function with H\"{o}lder function.
\begin{lemma}\label{multiplier estimate}
Let $\mathcal O\subset \mathbb R^d$ be an open bounded domain, and let
$0<s<\gamma<1$. Then if
$f\in C^{0,\gamma}(\overline{\mathcal O})$ and $g\in H^s(\mathcal O)$, then
$fg\in H^s(\mathcal O)$ and
$$
\|fg\|_{H^s(\mathcal O)}
\le c(s,\alpha,\mathcal{O})\|f\|_{C^{0,\gamma}(\overline{\mathcal O})} \|g\|_{H^s(\mathcal O)}.$$
\end{lemma}
\begin{proof}
Since $0<s<1$, 
$\|h\|_{H^s(\mathcal O)}^2 \simeq \|h\|_{L^2(\mathcal O)}^2+[h]_{H^s(\mathcal O)}^2
$, with
$$
[h]_{H^s(\mathcal O)}^2
=
\int_{\mathcal O}\int_{\mathcal O}
\frac{|h(x)-h(y)|^2}{|x-y|^{d+2s}}\d x\d y.
$$
It suffices to bound the seminorm $[fg]_{H^s(\mathcal{O})}$. In view of the inequality 
$|f(x)g(x)-f(y)g(y)|\le |f(x)||g(x)-g(y)|
+|f(x)-f(y)||g(y)|$, we obtain
$$
[fg]_{H^s(\mathcal O)}^2 \le 2{\rm I}+2{\rm II},
$$
with
$$
{\rm I}
=\int_{\mathcal O}\int_{\mathcal O}
\frac{|f(x)|^2|g(x)-g(y)|^2}{|x-y|^{d+2s}}
\d x\d y
\quad\mbox{and}\quad
{\rm II}
=\int_{\mathcal O}\int_{\mathcal O}
\frac{|f(x)-f(y)|^2|g(y)|^2}{|x-y|^{d+2s}}
\d x\d y.
$$
Since $f\in L^\infty(\mathcal O)$, we have
$$
{\rm I} \le\|f\|_{L^\infty(\mathcal O)}^2\int_{\mathcal O}\int_{\mathcal O} \frac{|g(x)-g(y)|^2}{|x-y|^{d+2s}} \d x\d y\le \|f\|_{L^\infty(\mathcal O)}^2 [g]_{H^s(\mathcal O)}^2 \le\|f\|_{C^{0,\gamma}(\overline{\mathcal O})}^2 \|g\|_{H^s(\mathcal O)}^2.
$$
The H\"older continuity of $f$ indicates
$|f(x)-f(y)|\le \|f\|_{C^{0,\gamma}(\overline{\mathcal O})}|x-y|^\gamma.$
Therefore,
\begin{align*}
{\rm II}
&\le
\|f\|_{C^{0,\gamma}(\overline{\mathcal O})}^2
\int_{\mathcal O}\int_{\mathcal O}
\frac{|x-y|^{2\gamma}|g(y)|^2}{|x-y|^{d+2s}}
\d x\d y\\
&=\|f\|_{C^{0,\gamma}(\overline{\mathcal O})}^2\int_{\mathcal O}
|g(y)|^2\left(\int_{\mathcal O}
|x-y|^{-d-2s+2\gamma}\d x\right)\d y.
\end{align*}
Since $\mathcal O$ is bounded and $\gamma>s$, there exists  $c=c(s,\gamma,\mathcal O)>0$ such that
$\int_{\mathcal O}
|x-y|^{-d-2s+2\gamma}\d x \le c$
uniformly for all $y\in\mathcal O$. Therefore,
$$
{\rm II} \le
c\|f\|_{C^{0,\gamma}(\overline{\mathcal O})}^2
\int_{\mathcal O}|g(y)|^2\,\d y\le
c\|f\|_{C^{0,\gamma}(\overline{\mathcal O})}^2
\|g\|_{L^2(\mathcal O)}^2.
$$
Combining the bounds on ${\rm I}$ and ${\rm II}$ gives the desired estimate.
\end{proof}

\section{Proof of Lemma \ref{lem:max-princ}}\label{app:subdiff}
First we state a nonnegativity result. The proof is essentially contained in \cite[Sections 2 and 3]{Zacher:2008}, but is given for the convenience of readers. Below the notation $(\cdot,\cdot)$ denotes the $L^2(\mathcal{O})$ inner product. 
\begin{lemma}\label{lem:comparison}
Let $u$ be the solution to problem \eqref{eq:tf-forward}. Under Assumption \ref{ass:subdiff}, for any $f\in\mathcal{F}$, there holds $u\ge0$ in $\mathcal{O}\times[0,T]$.
\end{lemma}
\begin{proof}
The proof follows that of \cite[Theorem 3.1]{Zacher:2008}, by rewriting the Caputo derivative $\partial_t^\alpha u$ in the form of the Riemann-Liouville derivative $\partial_t( \kappa\ast(u-u_0))$, with $\kappa(t)=g_{1-\alpha}(t)=\frac{t^{-\alpha}}{\Gamma(1-\alpha)}$ being the Riemann-Liouville integral kernel, and $\ast$ the convolution in time. Note that under Assumption \ref{ass:subdiff}, $u\in C([0,T];H^2(\Omega))$ \cite[Lemma 2.3(i)]{doi:10.1137/21M1446708}. Thus, we may equivalently write \eqref{eq:tf-forward} as
\begin{equation}\label{eqn:tf-forward:RL}
\partial_t\bigl(\kappa*(u-v)\bigr)-\Delta u + f\,u = h, \quad \text{in } \mathcal{O}\times (0,T].
\end{equation}
Let $\ell=g_\alpha$ be the resolvent kernel of $\kappa$, so $\kappa*\ell=1$ on $(0,\infty)$. For every $n\in\mathbb{N}$, we denote by $\kappa_n=\kappa*h_n\in W^{1,1}([0,T])$ the Yosida approximations of the kernel $\kappa$, with $h_n\ge0$, $\ell*\kappa_n=1*h_n$, and for any $1\leq p<\infty$,
\begin{equation}\label{eqn:approx-id} h_n*\varphi\to\varphi  \mbox{ in }L^p(0,T),\quad\forall \varphi\in L^p(0,T);
\end{equation}
see \cite[Section 2]{Zacher:2008} for further details. By \cite[Lemma 3.1]{Zacher:2008}, the weak formulation of \eqref{eqn:tf-forward:RL} is equivalent to
\begin{equation*}
\bigl(\partial_t(\kappa_n*(u-v)),\psi\bigr) + (h_n*\nabla u,\nabla\psi) + (h_n*(f u),\psi) = (h_n*h,\psi), \quad \forall \psi\in H_0^1(\mathcal{O}),
\end{equation*}
for a.e. $t\in(0,T)$, and all $n\in\mathbb{N}$. Now
fix $t$ and take $\psi=u^-:=-\max(-u,0)$. Since $u=g\ge0$ on $\partial\mathcal{O}$, we have $u^-=0$ on $\partial\mathcal{O}$. Hence,  $\psi\in H_0^1(\mathcal{O})$ is a valid test function, and
\begin{equation*}
\bigl(\partial_t(\kappa_n*(u-v)),-u^-\bigr) + (h_n*\nabla u,-\nabla u^-) + (h_n*(f u),-u^-) = (h_n*h,-u^-).
\end{equation*}
Since $\kappa_n\in W^{1,1}([0,T])$, using the convex inequality \cite[Corollary 4.1]{Zacher:2019} (with  $H(v)=(v^-)^2$) and the assumption $v\geq0$, we have
\[
-u^-\partial_t(\kappa_n*(u-v)) \ge \frac12\,\partial_t\bigl(\kappa_n*(u^-)^2\bigr).
\]
Meanwhile, since $h\geq 0$, $h_n\geq0$, $u^-\geq0$, we have $(h_n*h,\,u^-)\geq0$. Consequently,
\begin{equation*}
\frac12\,\partial_t\bigl(\kappa_n*(u^-)^2,\,1\bigr) - (h_n*\nabla u,\nabla u^-) - (h_n*(f u),u^-) \le 0.
\end{equation*} 
By convoluting the inequality with the kernel $\ell$, we obtain 
\begin{equation}\label{eqn:conv-ineq}
\frac12\,\bigl(\ell\ast(\partial_t(\kappa_n*(u^-)^2)),1\bigr) - \ell\ast(h_n*\nabla u,\nabla u^-) - \ell\ast(h_n*(f u),u^-) \le 0.
\end{equation}
Note that $k_n\ast (u^-)^2\in W^{1,1}([0,T];L^2(\Omega))$ with $(k_n\ast (u^-)^2)(0)=0$. Thus, from the relation $\ell*\kappa_n=h_n$, we deduce
\begin{equation*}
\ell\ast(\partial_t(\kappa_n\ast(u^-)^2)) = \partial_t(\ell\ast\kappa_n\ast(u^-)^2) =  \partial_t(h_n\ast(u^-)^2).
\end{equation*}
Then by integrating the inequality \eqref{eqn:conv-ineq} from $0$ to $t$ and noting $u^-(\cdot,0)=0$, we arrive at
\begin{equation}\label{eqn:weak-integrated}
\frac12\,\bigl(h_n*(u^-)^2,1\bigr)
- \int_0^t \ell*(h_n*\nabla u,\nabla u^-)\,{\rm d}s
- \int_0^t \ell*(h_n*(f u),u^-)\,{\rm d}s
\le 0.
\end{equation}
Note that the kernels $h_n$ form an approximate identity in time: for every $\varphi\in L^1(0,T;L^2(\mathcal{O}))$, $h_n*\varphi\to\varphi$ in $L^1(0,T;L^2(\mathcal{O}))$, cf. \eqref{eqn:approx-id}.
Using this property on the functions $\varphi=(u^-)^2$ and $\varphi=\nabla u$, $\varphi=f u$, we deduce for a.e. $s\in(0,t)$,
$(h_n*\nabla u,\,\nabla u^-) \to (\nabla u,\,\nabla u^-)$ and
$(h_n*(f u),\,u^-) \to (f u,\,u^-)$.
Since the convolution with $\ell\in L^1(0,T)$ is continuous on $L^1(0,T)$,
by passing to the limit in \eqref{eqn:weak-integrated}, we obtain
\begin{equation*}
\frac12\,\|u^-\|_{L^2(\mathcal{O})}^2
- \int_0^t \ell*(\nabla u,\nabla u^-)\,{\rm d}s
- \int_0^t \ell*(f u,u^-)\,{\rm d}s
\le 0.
\end{equation*}
On the support of $u^-$, $\nabla u = -\nabla u^-$ and $u = -u^-$. Hence,
\[
\frac12\|u^-\|_{L^2(\mathcal{O})}^2 + \int_0^t \ell*\|\nabla u^-\|_{L^2(\mathcal{O})}^2\,{\rm d}s
+ \int_0^t \ell*(f u^-,u^-)\,{\rm d}s
\le 0,
\]
which directly implies the estimate $\|u^-\|_{L^2(\mathcal{O})}^2 \le 0$. 
Thus $u^-=0$ or $u\geq 0$ a.e. in $\mathcal{O}\times(0,T]$. Under Assumption \ref{ass:subdiff}, we have $u\in C([0,T];H^2(\Omega))$ \cite[Lemma 2.3(i)]{doi:10.1137/21M1446708}, and by Sobolev embedding, $u\in C([0,T];C(\overline{\mathcal{O}}))$. Thus the desired assertion holds.
\end{proof}

Now we can prove Lemma \ref{lem:max-princ}.
\begin{proof}
Let $w:=\partial_t^\alpha u$. Then by \cite[(2.8)]{doi:10.1137/21M1446708}, $w$ solves 
\begin{equation}\label{eq:w-positive}
\left\{\begin{aligned}
\partial_t^\alpha w-\Delta w+fw&=0,\quad
    \text{in }\mathcal{O}\times(0,T),\\
w&=0,
    \quad \text{on }\partial\mathcal{O}\times(0,T),\\
w(\cdot,0)&=h +\Delta u_0 - f u_0,
    \quad \text{in }\mathcal{O}.
\end{aligned}\right.    
\end{equation}
Under Assumption \ref{ass:subdiff}, $w(\cdot,0)\in L^2(\mathcal{O})$. Then by the regularity theory for subdiffusion \cite[Lemma 2.3(i)]{doi:10.1137/21M1446708}, we have $w(t)\in C((0,T];H^2(\mathcal{O}))$, and thus by Sobolev embedding \cite{Adams:2003} (for $d=1,2,3$),
$$\partial_t ^\alpha u = w \in C((0,T];H^2(\mathcal{O})\cap H_0^1(\mathcal{O}))\subset C((0,T];C(\overline{\mathcal{O}})).$$
Meanwhile, by the definition of $\mathcal F$, $f\leq\Lambda$ a.e. in $\mathcal{O}$. This and Assumption \ref{ass:subdiff} yield 
$$w(\cdot,0) = h+\Delta u_0 - f u_0=(h+\Delta u_0-\Lambda u_0)+(\Lambda-f)u_0\geq 0.$$
Then by Lemma \ref{lem:comparison}, 
$w=\partial_t^\alpha u\geq 0$ a.e. in $\mathcal{O}\times(0,T)$.
Finally, using the elementary identity 
$$
u(x,t)=u(x,0)+\int_0^t\frac{(t-s)^{\alpha-1}}{\Gamma(\alpha)}\partial_s^\alpha u(x,s)\, {\rm d}s,
$$
we deduce $u(x,t)\geq u(x,0)\geq m_T$ for a.e. $(x,t)\in\mathcal{O}\times(0,T)$.
\end{proof}

\begin{remark}
In the literature, there are several versions of the maximum principle for the subdiffusion equation; see, e.g., \cite{Zacher:2008,LuchkoYamamoto:2019,jin2021fractional,Zacher:2019} and the references therein.  Lemma \ref{lem:comparison} is essentially identical with \cite[Theorem 3.2]{Zacher:2008}. Compared with \cite[Theorem 2.1]{LuchkoYamamoto:2017}, which requires $f\in \mathcal{A}\cap C(\overline{\mathcal{O}})$, Lemma \ref{lem:comparison} requires only $f\in \mathcal{A}$. 
\cite[Theorem 2.1]{LuchkoYamamoto:2017} was used in \cite{doi:10.1137/21M1446708,kow2025consistency} to obtain the relevant positivity result.
\end{remark}

\bibliographystyle{abbrv}
\bibliography{ref}

\end{document}